\documentclass[a4paper, 10pt,oneside]{amsart}
\usepackage{amsmath, amssymb, amsfonts}
\usepackage[all]{xy}
\usepackage{verbatim}
\newtheorem{tm}{Theorem}
\newtheorem{lm}[tm]{Lemma}
\newtheorem{prop}[tm]{Proposition}
\newtheorem{kor}[tm]{Corollary}

\newcommand{\st}    {\,|\,}

\newtheorem{napomena}[tm]{Remark}
\newenvironment{dokaz}
{\noindent\emph{Proof:}\ }
{\hfill $\square$}

\newcommand{\Z}
{{\mathbb Z}}
\newcommand{\N}
{{\mathbb N}}
\newcommand{\C}
{{\mathbb C}}

\newcommand{\baza}
{{\mathcal{B}}}
\newcommand{\g}
{{\mathfrak g}}

\newcommand{\gt}
{\tilde{{\mathfrak g}}}

\newcommand{\he}
{{\mathfrak h}^e}
\newcommand{\nt}
{\tilde{{\mathfrak n}}}
\newcommand{\h}
{{\mathfrak h}}
\newcommand{\n}
{{\mathfrak n}}
\newcommand{\gsl}
{{\mathfrak sl}}

\newcommand{\Gamt}
{\tilde{{\Gamma}}}

\begin{document}

\author{Goran Trup\v{c}evi\'{c}}
\title{Characters of Feigin-Stoyanovsky's type subspaces of level one modules for affine Lie algebras of types $A_\ell^{(1)}$ and $D_4^{(1)}$}
\address{Department of Mathematics, University of Zagreb, Bijeni\v cka 30, Zagreb, Croatia}
\curraddr{}
\email{gtrup@math.hr}
\thanks{}
\subjclass[2000]{Primary 17B67; Secondary 05A19.\\ \indent Partially supported by the Ministry of Science and
Technology of the Republic of Croatia, Project ID 037-0372794-2806}
\keywords{affine Lie algebras, principal subspaces, character formulas}
\date{}
\dedicatory{}

\begin{abstract}
We use combinatorial description of bases of Feigin-Stoyanovsky's
type subspaces of standard modules of level $1$ for affine Lie
algebras of types $A_\ell^{(1)}$ and $D_4^{(1)}$ to obtain character
formulas. These descriptions naturally lead to systems of recurrence
relations for which we also find solutions.
\end{abstract}

\maketitle


\section{Introduction}

Principal subspaces were introduced  by B.L.\,Feigin and
A.\,Stoyanovsky in [FS] where they gave a construction of bases of
standard modules $L(\Lambda)$ consisting of semi-infinite monomials
and monomial bases of their principal subspaces, and also calculated
characters of both principal subspaces and the whole standard
modules for affine Lie algebra $\gt$ of  type $A_1^{(1)}$. A similar
approach was used by M.\,Primc  in [P1,2] where he constructed
semi-infinite monomial bases for all standard modules for affine Lie
algebras of  type $A_\ell^{(1)}$ and for basic modules
$L(\Lambda_0)$ for any classical affine Lie algebra. Instead of
principal subspaces of Feigin and Stoyanovsky, Primc used so-called
Feigin-Stoyanovsky's type subspace. Later, in [FJLMM] it was noted
that bases of Feigin-Stoyanovsky's type subspaces from [P1] were
parameterized by $(k,\ell+1)$-admissible configurations which were
studied in [FJLMM], [FJMMT1,2].

G.\, Georgiev generalized Feigin-Stoyanovsky's results to a certain
class of standard modules for affine Lie algebras of type
$A_\ell^{(1)}$ (see [G]). In the proof of linear independence, Georgiev used
intertwining operators between standard modules. S.\,Capparelli,
J.\,Lepowsky and A.\,Milas in [CLM1,2] used intertwining operators
to obtain exact sequences of principal subspaces and recurrence
relations for their characters. This approach was further
investigated in [C1,2] and [CalLM1-3].

Motivated by Georgiev's and Capparelli-Lepowsky-Milas' way of using
intertwining operators, Primc in [P3]  gave a simpler proof of
linear independence of bases from [P1], and in [T1,2] and [B1,2] new
constructions of bases in $A_\ell^{(1)}$ and $D_\ell^{(1)}$ cases
were given. Furthermore, M.\,Jerkovi\'c in [J1] used the proof of
linear independence from [P3] to obtain exact sequences of
Feigin-Stoyanovsky's type subspaces and recurrence relations for the
corresponding characters in the  $A_\ell^{(1)}$-case. By solving
these relations, Jerkovi\'c in [J2] obtained character formulas in
the $A_2^{(1)}$-case, which agreed with formulas from [FJMMT1,2].

In this paper we use combinatorial description of bases of
Feigin-Stoyanovsky's type subspaces of standard modules of level $1$
from [P2], [T1] and [B1] to obtain character formulas. These
descriptions naturally lead to systems of recurrence relations for
which we also find solutions.

Let $\g$ be a simple complex Lie algebra, $\h\subset\g$ its Cartan
subalgebra, $R$ the corresponding root system. Let $\g=\h+\sum_{\alpha\in R}\g_\alpha$
be a root decomposition of $\g$. Fix root vectors
$x_\alpha\in\g_\alpha$. Let $\langle
\cdot,\cdot\rangle$ be a normalized invariant bilinear form on
$\g$, and by the same symbol denote the induced form on $\g^*$.
Let $\Pi=\{\alpha_1,\dots,\alpha_\ell\}$ be a basis of the root system $R$, and $\{\omega_1,\dots,\omega_\ell\}$ the
corresponding set of fundamental weights. Fix a minuscule fundamental weight
$\omega$ and set $\Gamma=\{\gamma\in
R\,|\,\langle\gamma,\omega\rangle=1\}$, $\g_1=\sum_{\alpha \in \Gamma}\, \g_\alpha$. The set $\Gamma$ is called
{\em the set of  colors}.

Let $\gt=\g\otimes
\C[t,t^{-1}]\oplus \C c \oplus \C d$ be the associated affine Lie algebra, where $c$ is the canonical
central element, and $d$ is the degree operator. Elements
$x_\alpha(r)=x_\alpha\otimes t^r$ are fixed real root vectors.
Let $\gt_1=\g_1\otimes\C[t,t^{-1}]$, a commutative Lie
subalgebra with a basis
$\{x_\gamma(-r)\,|\,r\in\Z,\gamma\in\Gamma\}$.
%
%
Let $L(\Lambda)$ be a standard $\gt$-module of level $1$, with a
fixed highest weight vector $v_\Lambda$. A {\em Feigin-Stoyanovsky's
type subspace} of $L(\Lambda)$ is a $\gt_1$-submodule of
$L(\Lambda)$ generated with $v_\Lambda$,
$$W(\Lambda)=U(\gt_1)\cdot v_\Lambda\subset L(\Lambda).$$

For the Lie algebra $\g$ of  type $A_\ell$
it was shown in [P2] and [T1] that monomial vectors $\underline{x} v_\Lambda$, where $\underline{x}=x_{\gamma_n}(-r_n)
\cdots  x_{\gamma_1}(-r_1)$, $\gamma_i\in\Gamma$, $r_i\in\N$, such that $\underline{x}$ satisfy
certain combinatorial conditions called {\it difference} and {\it initial conditions}, constitute
a basis of $W(\Lambda)$. The analogous fact was proved in [P2] and [B1] for $\g$
of  type $D_\ell$.

To obtain character formula when Lie algebra $\g$ is of  type $A_\ell$, we first consider two particular cases,
when $\omega=\omega_1$ and $\omega=\omega_\ell$; these are the cases that were considered in
[P1], [FJLMM], [FJMMT1,2] and [J1-3], but for higher-level modules. 
For every $\h$-weight subspace of $W(\Lambda)$, we construct a bijection between
the basis of that subspace and products of partitions of certain length.
This gives formulas \eqref{kar1r_jed} and \eqref{kar2r_jed} for the character of $W(\Lambda)$, that were
already known (e.g. in [J3]).
In the case $\omega=\omega_m$, $1<m<\ell$, the set of colors $\Gamma$ can be decomposed into a product
of "rows" and "columns". The sets of rows and columns can be regarded as sets of colors for the two
particular cases that have already been  considered. For a given $\h$-weight subspace of $W(\Lambda)$,
we consider its basis elements $\underline{x} v_\Lambda$, where $\underline{x}=x_{\gamma_n}(-r_n)
\cdots  x_{\gamma_1}(-r_1)$, $\gamma_i\in\Gamma$, $r_i\in\N$. To every such basis element we can attach
its path ${\bf p}(\underline{x})=(\gamma_n,\dots,\gamma_1)$, and conversely, to every path $\bf p$
we can attach a basis element $\underline{x}({\bf p})$ that will be minimal in some sense. By the decomposition
of $\Gamma$, for every path $\bf p$ in $\Gamma$,
we have the corresponding paths of rows and columns. We use character formulas for the two particular cases
to find ``graded cardinality'' of the set of ``minimal'' monomials for paths corresponding to the given
$\h$-weight. From this we obtain formula \eqref{KarAm_jed} for the character of $W(\Lambda)$.

When Lie algebra $\g$ is of  type $D_4$, we decompose the set of colors into two subsets
that correspond to the cases $A_2$, with $\omega=\omega_2$, and $A_3$, with $\omega=\omega_2$.
We use
character formulas \eqref{kar1r_jed} and \eqref {KarAm_jed} for the latter cases to obtain character formula \eqref{KarD4_jed} in the $D_4$-case.

Both in $A_\ell$ and $D_4$ cases, descriptions of combinatorial
bases naturally lead to systems of recurrence relations. We can find
solutions of these systems in a similar way to the one we used for
calculating character formulas of Feigin-Stoyanovsky's type
subspaces.

The outline of this paper is as follows: in Section \ref{FS_sect} we
introduce basic definitions. In Section \ref{al_sect} we find
character formulas in the $A_\ell$-case. We also find solutions of
the corresponding system of recurrence relations. In Section
\ref{d4_sect} we do the same thing in the $D_4$-case.

\section{Feigin-Stoyanovsky's type subspace}

\label{FS_sect}

Let $\g$ be a simple finite-dimensional Lie algebra. Let $\h\subset\g$ be a Cartan subalgebra of
$\g$ and $R$ the corresponding root system. Fix a basis
$\Pi=\{\alpha_1,\dots,\alpha_\ell\}$ of $R$. Then we have the
root decomposition $\g=\h \oplus \coprod_{\alpha\in R} \g_\alpha$ and the triangular decomposition $\g=\n_-\oplus \h \oplus \n_+$. Let
$\theta=k_1\alpha_1+\dots+k_\ell \alpha_\ell$ be the maximal root. Let  $\langle
\cdot,\cdot\rangle$ be a a normalized invariant bilinear form on
$\g$ such that $\langle
\theta,\theta \rangle=2$; we identify $\h$ with $\h^*$ via $\langle
\cdot,\cdot\rangle$. For $\alpha\in R$ let $\alpha^\vee=2\alpha/\langle \alpha,\alpha\rangle$ denote
the corresponding coroot. Also for each root $\alpha\in R$ fix a root vector $x_\alpha\in\g_\alpha$.
%
Let $\{\omega_1,\dots,\omega_\ell\}$ be the set of
fundamental weights of $\g$,
$\langle\omega_i,\alpha_j\rangle=\delta_{ij},\,i,j=1,\dots,\ell$.
Denote by $Q=\sum_{i=1}^\ell\Z \alpha_i$
the root lattice, and by $P=\sum_{i=1}^\ell\Z \omega_i$ the weight lattice
 of $\g$. Denote by $P^+=\sum_{i=1}^\ell\Z_{\geq 0} \omega_i$ the set of dominant integral weights.

Let $\gt$ be the associated untwisted affine Lie algebra,
$$\gt=\g\otimes \C[t,t^{-1}]\oplus \C c \oplus \C d,$$
with commutation relations
\begin{eqnarray*}
 & & \hspace{-3ex} [x(i),y(j)]  =  [x,y](i+j)+ i\langle x,y\rangle \delta_{i+j,0}c, \\
 & & \hspace{-3ex} [c,\gt]  =  0,\quad [d,x(j)]  =  j x(j),
\end{eqnarray*}
where $x(j)=x\otimes t^j$ for $x\in\g$, $j\in\Z$ (cf. [K]).

Set $\he=\h\oplus\C c \oplus\C d,\, \nt_\pm=\g\otimes t^{\pm 1}\C
[t^{\pm 1}]\oplus \n_\pm$.  Then $\gt$ also has the triangular decomposition
$\gt=\nt_-\oplus \he \oplus \nt_+$.
Usual extensions of bilinear
forms $\langle\cdot,\cdot\rangle$ onto $\he$ and $(\he)^*$ are
denoted by the same symbols (we take $\langle c,d \rangle=1$). 
Denote by $ \alpha_0, \alpha_1, \dots, \alpha_\ell \in
(\he)^{*}$ the simple roots, and by $\Lambda_0, \Lambda_1, \dots,
\Lambda_\ell \in (\he)^{*}$ the corresponding fundamental weights.  Then
$\Lambda_0(c) =1$, $\Lambda_i(c) =k_i$ for $i=1, \dots, \ell$.

%

Weight $\omega\in P$ is said to be \emph{minuscule} if
$\langle\omega,\alpha\rangle\in \{-1,0,1\}$ for $\alpha \in R$.
A dominant integral weight $\omega\in P^+$
is minuscule if and only if
$\langle\omega,\theta\rangle = 1.$

Fix a minuscule weight $\omega\in P$.
Set
$$\Gamma =
\{\,\alpha \in R \mid \langle\alpha,\omega\rangle = 1\}.$$
Then
$$
\mathfrak g  =
\mathfrak g_{-1} \oplus \mathfrak g_0 \oplus \mathfrak g_1, $$
where
$$ {\mathfrak g}_0  =   {\mathfrak h} \oplus
\sum_{\langle \alpha, \omega \rangle =0}\, {\mathfrak g}_\alpha,\qquad
 {\mathfrak g}_{\pm1} =
\sum_{\alpha \in \pm \Gamma}\, {\mathfrak g}_\alpha,
$$ is a $\mathbb Z$-gradation of ${\mathfrak g}$.
Subalgebras ${\mathfrak g}_1$ and ${\mathfrak g}_{-1}$ are
commutative.
We call elements  $\gamma\in\Gamma$ {\em
colors} and the set $\Gamma$ {\em the  set of colors}.

 The $\mathbb Z$-gradation of ${\mathfrak g}$ induces the $\Z$-gradation of  affine Lie
algebra $\gt$
\begin{eqnarray*}
 & \gt = \gt_{-1} + \gt_0 + \gt_1, &
\\
 & \gt_0 = {\mathfrak g}_0\otimes\C [t,t^{-1}]\oplus \C c \oplus \C d,\qquad
\gt_{\pm 1} = {\mathfrak g}_{\pm 1}\otimes\C [t,t^{-1}]. &
\end{eqnarray*}
Again, $\gt_{-1}$ and $\gt_1$ are commutative subalgebras. Set $\gt_1^-=\gt_1\cap \nt_-$.

Let $L(\Lambda_k)$ be a standard (i.e. integrable highest weight) $\gt$-module
of level $\Lambda_k (c)=1$. Denote by $v_{\Lambda_k}$ the highest weight vector of $L(\Lambda_k)$.
Define a
\emph{Feigin-Stoyanovsky's type subspace}
$$W(\Lambda_k)=U(\gt_1)\cdot v_{\Lambda_k}=U(\gt_1^-)\cdot v_{\Lambda_k}\subset L(\Lambda_k).$$

Set
$$\Gamt=\{x_\gamma(-r) \mid \gamma\in\Gamma,r\in \Z\},\quad \Gamt^-=\{x_\gamma(-r) \mid \gamma\in\Gamma,r\in \N\}.$$
Since the subalgebra $\gt_1$ is commutative, we have $U(\gt_1)\cong \C[\Gamt]$ and $U(\gt_1^-)\cong \C[\Gamt^-]$. We often refer to
elements of $\Gamt$ as to {\em variables}, {\em elements} or {\em factors}
of a monomial from $U(\gt_1)$.

\section{The case $A_\ell$, $\ell\geq 1$}
\label{al_sect}

Let ${\mathfrak g}$ be a simple finite-dimensional Lie algebra of  type $A_\ell$. In this case all fundamental weights are minuscule.
Fix a minuscule weight $\omega=\omega_m$, $m\in\{1,\dots,\ell\}$.
The set of colors $\Gamma$  is parameterized by two sets of
indices
$$ \Gamma=\{(ij)\,|\, i=1,\ldots,m; j=m,\ldots,\ell\}, $$
where
\begin{equation}
\label{gij_jed}  (ij)=\alpha_i+\cdots+\alpha_m+\cdots+\alpha_j,
\end{equation}
and thus we can think of it as a rectangle with rows ranging from $1$ to
$m$, and columns ranging from $m$ to $\ell$ (see Figure 1 in [T1]).
By $x_{ij}\in\g$ we denote the fixed root vector corresponding to the color
$(ij)$.

Linear order $<$ on the set of colors $\Gamma$ is defined as follows:
$(ij)<(i'j')$ if either
 $i>i'$ or  $i=i'$ and $j<j'$.
On the set of variables
$\Gamt$ we
define a linear order by:
$x_\gamma(-r)<x_{\gamma'}(-r')$ if
either $-r<-r'$ or $r=r'$ and $\gamma<\gamma'$.
Since the algebra $\gt_1$ is commutative, we assume that variables in monomials from $\C[\Gamt]$  are sorted ascendingly from left to right.

Let $L(\Lambda_k)$, $k=0,\dots,\ell$ be a standard $\gt$-module
of level $1$.
We use a description of a combinatorial basis of $W=W(\Lambda_k)$ from [P2] and [T1].
Define an energy function $E:\Gamma\times\Gamma\to\{0,1,2\}$ by
\begin{equation} \label{Energy_jed}
E((i'j'),(ij))=\left\{\begin{array}{ll}
0,  &  i'>i, j'<j,\\
1,  &  i'\leq i, j'<j \quad\textrm{or}\quad i'>i,j'\geq j,\\
2,  &  i'\leq i, j'\geq j.
\end{array}\right.
\end{equation}
Define $\theta:\Z\to \{0,1\}$ by $$\theta(n)=\left\{\begin{array}{ll}
0, & n<0,  \\
1, & n\geq 0.
\end{array}\right.$$
Then
\begin{equation}
\label{EnerAProd_jed}
E((i'j'),(ij))=\theta(i-i')+\theta(j'-j).
\end{equation}

We say that a monomial $$\underline{x}=x_{i_n j_n}(-r_n)\cdots x_{i_1 j_1}(-r_1)\in\C[\Gamt^-]$$ satisfies \emph{difference conditions}, or DC for short, if
$$r_{t+1}-r_t \geq E((i_{t+1}j_{t+1}),(i_t j_t)).$$
We say that $\underline{x}$ satisfies
{\em initial conditions} for $L(\Lambda_k)$, or IC for short, if either $r_1\geq 2$ or
$r_1=1$ and either $i_1>k$, for $1\leq k\leq m$, or $j_1<k$, for $m \leq k \leq \ell$. Define
\begin{equation}\label{baza_def}
\baza_W=\{\underline{x}\in\C[\Gamt^-] \,|\, \underline{x} \ \textrm{satisfies DC and IC for}\ L(\Lambda_k)\}.
\end{equation}

\begin{tm}
The set $\{\underline{x} v_{\Lambda_k} \st \underline{x}\in\baza_W\}$
is a basis of $W$.
\end{tm}

For a monomial $\underline{x}=x_{\gamma_n}(-r_n)\cdots x_{\gamma_1}(-r_1) \in \C[\Gamt^-]$, define {\em weight} and {\em degree} by
$$w(\underline{x})=\gamma_1+\dots+\gamma_n,\qquad d(\underline{x})=r_1+\dots+r_n.$$
For $\alpha\in P$, set $\underline{z}^\alpha=z_1^{\langle\alpha,\omega_1\rangle}\cdots z_\ell^{\langle\alpha,\omega_\ell\rangle}$.
The {\em character} of $W$ is the formal sum
$$\chi_{W}(z_1,\dots,z_\ell,q)=\sum_{\underline{x}\in\baza_W} q^{d(\underline{x})} \underline{z}^{w(\underline{x})}.$$
For a fixed $\alpha=n_ 1 \alpha_1+\dots+n_\ell \alpha_\ell\in P^+$, define
$\baza^\alpha_W  = \{ \underline{x} \in\baza_W \st  w(\underline{x})=\alpha \}$
and $\chi_{W}^\alpha(q)=\sum_{\underline{x}\in \baza^\alpha_W} q^{d(\underline{x})}$.
Obviously, $\chi_{W}(z_1,\dots,z_\ell,q)=\sum_{\alpha\in P^+} \chi_{W}^\alpha(q)
\underline{z}^\alpha$.
We sometimes use symbols $\baza_W^{n_1,\dots,n_\ell}$, $\chi_W^{n_1,\dots,n_\ell}(q)$
instead of $\baza_W^\alpha$ and $\chi_W^\alpha(q)$.

From $\eqref{gij_jed}$ it immediately follows that $\chi_W^\alpha(q)=0$ unless $0\leq n_1\leq\dots\leq n_m\geq\dots\geq n_\ell\geq 0$.

A nondecreasing sequence of nonnegative integers $\lambda = (\lambda_1,\dots,\lambda_n)$, $0\leq\lambda_1\leq\dots\leq\lambda_n$ is called
{\em a partition} of length at most $n$. The  sum $|\lambda|=\sum_i \lambda_i$ is called  {\em weight}  of $\lambda$. Denote by $\pi_n$ the
set of partitions of length at most $n$.

For a monomial $\underline{x}=x_{\gamma_n}(-r_n)\cdots x_{\gamma_1}(-r_1) \in \C[\Gamt^-]$ and a partition $\lambda\in\pi_n$ define monomials
\begin{eqnarray}
\nonumber \underline{x}^{\pm\triangledown} \hspace{-1.7ex} & = & \hspace{-1.7ex} x_{\gamma_n}(-r_n \pm (n-1))\cdots x_{\gamma_2}(-r_2 \pm 1)x_{\gamma_1}(-r_1),\\
\nonumber \underline{x}^{\pm} \hspace{-1.7ex} & = & \hspace{-1.7ex} x_{\gamma_n}(-r_n \pm 1)\cdots x_{\gamma_1}(-r_1 \pm 1), \\
\nonumber \underline{x}^{\pm r} \hspace{-1.7ex} & = & \hspace{-1.7ex} x_{\gamma_n}(-r_n \pm r)\cdots x_{\gamma_1}(-r_1 \pm r),\quad \textrm{for}\ r\in\N, \\
\label{ParMon_jed} \underline{x}(\lambda) \hspace{-1.7ex} & = & \hspace{-1.7ex} x_{\gamma_n}(-n-\lambda_{n}) \cdots
x_{\gamma_1}(-1-\lambda_{1}).
\end{eqnarray}
We emphasize that the monomial $\underline{x}$ is assumed to be sorted ascendingly from left to right.
Note that if $\underline{x}$ satisfies difference and initial conditions,
then the variables in $\underline{x}$ are sorted ascendingly from left to right.

 \subsection{Character formula in the case $\omega=\omega_1$ or $\omega=\omega_\ell$}
\label{Al1_sect}

Consider the second case, $\omega=\omega_\ell$; the first case can be treated analogously. Fix $W=W(\Lambda_k)$, $0\leq k\leq\ell$.

The set of colors in this case is  $\Gamma=\{(1 \ell),\dots, (\ell \ell)\}$. For simplicity, we write $(i)$ and $x_i$ instead of $(i \ell)$ and
$x_{i \ell}$, for $i=1,\dots,\ell.$
The formula \eqref{Energy_jed} for the energy function in this case takes a simpler form:
\begin{equation} \label{EnerAl_jed}
E((i'),(i))  =
\theta(i-i')+1.
\end{equation}
Set $E'((i'),(i))=E((i'),(i))-1=\theta(i-i')$. We say that a monomial $\underline{x}=x_{i_n}(-r_n)\cdots x_{i_1}(-r_1)\in\C[\Gamt^-]$ satisfies DC' if
$r_{t+1}-r_t \geq E'((i_{t+1}),(i_t)).$
The following lemma is obvious
\begin{lm} \label{DCprime_lm}
A monomial $\underline{x}$ satisfies DC if and only if  $\underline{x}^{+\triangledown}$ satisfies DC'.
\end{lm}

Fix $0\leq n_1\leq n_2\leq \dots \leq n_\ell$ and set $\alpha=n_ 1 \alpha_1+\dots+n_\ell \alpha_\ell$. Set $n_i'=n_i-n_{i-1}$,
for $i=2,\dots,\ell$, and $n_1'=n_1$; then $\alpha=n_ 1' (1)+\dots+n_\ell' (\ell)$.
Let $\underline{\lambda}=(\lambda^1,\dots,\lambda^\ell)\in
\pi_{n_1'}\times\cdots\times\pi_{n_\ell'}$. For $i=1,\dots,\ell$, define
$\underline{x}_i=x_i(-n_i'-\theta(k-i)) x_i(-n_i'+1-\theta(k-i))\cdots
x_i(-1-\theta(k-i))$. Set $\underline{x}(\underline{\lambda})=\underline{x}_1(\lambda^1) \cdots \underline{x}_\ell (\lambda^\ell)$
, and reorder variables so that they are sorted ascendingly from left to right. Then obviously $\underline{x}(\underline{\lambda})$
satisfies DC' and IC for $L(\Lambda_k)$. Hence, by lemma \ref{DCprime_lm}, $\underline{x}(\underline{\lambda})^{-\triangledown}\in \baza_W^\alpha$.

Conversely, let  $\underline{x}\in \baza_W^\alpha$. Set $\underline{x}'=\underline{x}^{+\triangledown}$. Factorize
$\underline{x}'=\underline{x}_1 \cdots \underline{x}_\ell$ so that
$\underline{x}_i=x_i(-r_{n_i'}^i)\cdots x_i(-r_{1}^i)$ and $r_{n_i'}^i>\dots>r_{1}^i>0$ (this is possible since $E'((i),(i))=1$).
Define $\lambda_t^i=r_t^i-t$ for $i=1,\dots,\ell$, $t=1,\dots,n_i'$. Then obviously $\lambda^i=(\lambda_1,\dots,\lambda_{n_i'})\in\pi_{n_i'}$.
We have proved

\begin{tm}
The map
\begin{eqnarray*}
\pi_{n_1}\times \pi_{n_2-n_1}\cdots\times\pi_{n_\ell-n_{\ell-1}} \hspace{-1.7ex} & \to & \hspace{-1.7ex} \baza_W^\alpha,\\
\underline{\lambda}  \hspace{-1.7ex} & \mapsto & \hspace{-1.7ex} \underline{x}(\underline{\lambda})^{-\triangledown}
\end{eqnarray*}
is a bijection.
\end{tm}

Obviously
\begin{eqnarray*}
d(\underline{x}(\underline{\lambda})^{-\triangledown}) \hspace{-1.7ex} & = &
\hspace{-1.7ex} |\lambda| + \sum_{i=1}^\ell \frac{n_i'(n_i'+1)}{2}+\frac{(n_1'+\dots+n_\ell')(n_1'+\dots+n_\ell'-1)}{2}+\sum_{i=1}^k n_i'\\
 \hspace{-1.7ex} & = & \hspace{-1.7ex} |\lambda| + \sum_{i=1}^\ell n_i^2-\sum_{i=1}^{\ell-1} n_i n_{i+1}+n_k.
\end{eqnarray*}
As a consequence, we have

\begin{kor} \label{KarAl_kor}
For $0\leq n_1\leq n_2\leq \dots \leq n_\ell$, $\alpha=n_ 1 \alpha_1+\dots+n_\ell \alpha_\ell$,
\begin{equation} \label{kar1r_jed}
\chi_{W(\Lambda_k)}^\alpha(q)=\frac{ q^{\sum_{i=1}^\ell n_i^2-\sum_{i=1}^{\ell-1} n_i n_{i+1}+n_k}} {(q)_{n_1}(q)_{n_2-n_1}\cdots(q)_{n_\ell-n_{\ell-1}}},
\end{equation}
where $(q)_n=(1-q)\cdots(1-q^n)$.
\end{kor}

\begin{napomena} \label{KarA1_nap}
Analogous formula can be obtained in the $\omega=\omega_1$ case; for $n_1\geq n_2\geq \dots \geq n_\ell \geq 0$ we have
\begin{equation}
\label{kar2r_jed}
\chi_{W(\Lambda_k)}^{\alpha}(q)=\frac{ q^{\sum_{i=1}^\ell n_i^2-\sum_{i=1}^{\ell-1} n_i n_{i+1}+n_k}} {(q)_{n_\ell}(q)_{n_{\ell-1}-n_\ell}\cdots(q)_{n_1-n_2}},
\end{equation}
\end{napomena}

\subsection{Character formula in the case $\omega=\omega_m,\,1<m<\ell$}
\label{Al2_sect}

Define Lie subalgebras 
 $\g'=\langle x_{\alpha_i}, x_{-\alpha_i}, \alpha_i \st  i=1,\dots,m  \rangle$
 and $\g''=\langle x_{\alpha_i}, x_{-\alpha_i}, \alpha_i \st i=m,\dots,\ell  \rangle$  of types $A_m$ and $A_{\ell-m+1}$, respectively.
 We regard $\alpha_1,\dots,\alpha_m$ and $\alpha_m,\dots,\alpha_\ell$ as root bases, and $\omega_1,\dots,\omega_m$ and $\omega_m,\dots,\omega_\ell$ as fundamental
 weights for these subalgebras. Also, we regard $\Lambda_0,\Lambda_1,\dots,\Lambda_m$ and $\Lambda_0,\Lambda_m,\dots,\Lambda_\ell$ as fundamental weights for the
 corresponding affine Lie algebras. It will be clear from the context when the symbols $L(\Lambda_k)$, $W(\Lambda_k)$, $k=0,\dots,\ell$ denote the standard
module and the corresponding Feigin-Stoyanovsky'y type subspace for
$\g$, when for $\g'$ and when for $\g''$.

The set of colors $\Gamma$ is parameterized by two sets of indices -- the set of row-indices $\Gamma_1=\{1,\dots,m\}$ and the set of column-indices
$\Gamma_2=\{m,\dots,\ell\}$. We regard these two sets as sets of colors for $\g'$ and $\g''$, for the choice of minuscule weight $
\omega=\omega_m$ in both cases. Energy functions for $\Gamma_1$ and $\Gamma_2$ are
\begin{equation}\label{EnergyDecomp_jed}
E_1(i',i)=\theta(i-i')+1,\quad E_2(j',j)=\theta(j'-j)+1,
\end{equation}
for $i,i'\in\Gamma_1$ and $j,j'\in\Gamma_2$ (see \eqref{EnerAl_jed}).
By \eqref{EnerAProd_jed}, we have
$$E((i'j'),(ij))=E_1(i',i)+E_2(j',j)-2.$$

We consider the case $W=W(\Lambda_0)$ in detail, the other cases work in the analogous manner.

A {\em path} is a finite sequence of colors ${\bf p}=(\gamma_n,\dots,\gamma_1)$.
The number $l({\bf p})=n$ is called {\em length} of ${\bf p}$. The sum $w({\bf p})=\gamma_1+\dots+\gamma_n$ is called {\em weight} of ${\bf p}$.

To each monomial $\underline{x}=x_{\gamma_n}(-r_n)\cdots x_{\gamma_1}(-r_1) \in \C[\Gamt^-]$ we attach its path ${\bf p}(\underline{x})=(\gamma_n,\dots,\gamma_1)$.
Obviously $w({\bf p}(\underline{x}))=w(\underline{x})$.

Conversely, to a fixed path ${\bf p}=(\gamma_n,\dots,\gamma_1)$ we attach a monomial
$\underline{x}({\bf p})=x_{\gamma_n}(-r_n)\cdots x_{\gamma_1}(-r_1)$ such that
\begin{equation}\label{MinMon_jed}
r_1=1,\quad r_t=r_{t-1}+E(\gamma_t,\gamma_{t-1})\ \textrm{for}\ t=2,\dots, n.
\end{equation}
This is the ``minimal'' monomial of path $\bf p$ that satisfies difference and initial conditions for $L(\Lambda_0)$.
By this we mean that if $\lambda=(\lambda_1,\dots,\lambda_n)\in\pi_n$ is a partition of length at most $n$, then the
monomial $(\underline{x}({\bf p}))(\lambda)=x_{\gamma_n}(-r_n-\lambda_n)\cdots x_{\gamma_1}(-r_1-\lambda_1)$ also satisfies
difference and initial conditions, and all monomials of path $\bf p$ that satisfy difference and initial conditions can be obtained in this way.

Fix $0\leq n_1\leq\dots\leq n_m\geq\dots\geq n_\ell\geq 0$ and set $\alpha=n_ 1 \alpha_1+\dots+n_\ell \alpha_\ell$. The argument
from the preceding paragraph shows that
\begin{equation} \label{KarMin_jed}
\chi_{W(\Lambda_0)}^\alpha(q)=\frac{1}{(q)_{n_m}} \sum_{{\bf p},\,w({\bf p})=\alpha}q^{d(\underline{x}({\bf p}))},
\end{equation}
since $l({\bf p})=n_m$ for a path $\bf p$ of weight $\alpha$.

Fix a path ${\bf p}=((i_{n_m} j_{n_m}),\dots,(i_1 j_1))$ in $\Gamma$ of weight $\alpha$. Then, by \eqref{gij_jed}
\begin{eqnarray}
\label{retci_jed} & & \hspace{-3ex}
n_1=\#\{t|i_t=1\},\quad n_s-n_{s-1}=\#\{t|i_t=s\},\ \textrm{for}\ s=2,\dots,m,
 \\
\label{stupci_jed} & & \hspace{-3ex}
n_\ell=\#\{t|j_t=\ell\},\quad  n_s-n_{s+1}=\#\{t|j_t=s\},\ \textrm{for}\ s=m,\dots,\ell-1.
\end{eqnarray}
Denote by
${\bf p}_1=(i_{n_m},\dots,i_1)$ and ${\bf p}_2=(j_{n_m},\dots,j_1)$ the corresponding paths in $\Gamma_1$ and $\Gamma_2$.
Weights of ${\bf p}_1$ and ${\bf p}_2$ are $\alpha'=n_ 1 \alpha_1+\dots+n_m \alpha_m$ and $\alpha''=n_m \alpha_m+\dots+n_\ell \alpha_\ell$.
Conversely, if ${\bf p}_1$ and ${\bf p}_2$ are paths in $\Gamma_1$ and $\Gamma_2$ of weights $\alpha'$ and $\alpha''$, respectively, then
the corresponding path in $\Gamma$ will be of weight $\alpha$ (cf. \eqref{gij_jed}, \eqref{retci_jed}, \eqref{stupci_jed}).

Let $\underline{x}({\bf p})=x_{i_{n_m} j_{n_m}}(-r_{n_m})\cdots x_{i_{1} j_{1}}(-r_{1})$,
$\underline{x}({\bf p}_1)=x_{i_{n_m} }(-r_{n_m}')\cdots x_{i_{1} }(-r_{1}')$,
$\underline{x}({\bf p}_2)=x_{j_{n_m}}(-r_{n_m}'')\cdots x_{j_{1}}(-r_{1}'')$
be like in \eqref{MinMon_jed}. Then, by
\eqref{EnerAProd_jed}, \eqref{EnergyDecomp_jed} and \eqref{MinMon_jed}, we have
\begin{equation}\label{MinMonPocDeg_jed}
r_1=r_1'=r_1''=1
\end{equation}
and
\begin{equation} \label{MinDegRek_jed}
\left\{
\begin{array}{rcl}
r_t \hspace{-1.5ex} & = & \hspace{-1.5ex} r_{t-1} +\theta(i_{t-1}-i_t) + \theta(j_t-j_{t-1}),\\
r_t' \hspace{-1.5ex} & = & \hspace{-1.5ex} r_{t-1}' +\theta(i_{t-1}-i_t) + 1,\\
r_t'' \hspace{-1.5ex} & = & \hspace{-1.5ex} r_{t-1}' + \theta(j_t-j_{t-1}) + 1,
\end{array}\right.\end{equation}
for $t=2,\dots,n_m$. By induction, from this we obtain
\begin{equation}\label{MinDeg_jed}
r_t=r_t'+r_t''-2t+1,\qquad \textrm{for}\ t=1,\dots,n_m.
\end{equation}
This implies
\begin{equation}
\label{DegProdA_jed}
d(\underline{x}({\bf p})) = d(\underline{x}({\bf p}_1))+d(\underline{x}({\bf p}_2))-n_m^2.
\end{equation}
Consequently
\begin{eqnarray}
\label{KarMin2_jed}
 \sum_{{\bf p},\,w({\bf p})=\alpha}q^{d(\underline{x}({\bf p}))}
 \hspace{-1.7ex} & = & \hspace{-1.7ex}  \sum_{\substack{
{\bf p}_1,\,w({\bf p}_1)=\alpha'\\
{\bf p}_2,\,w({\bf p}_2)=\alpha''}} q^{d(\underline{x}({\bf p}_1))+d(\underline{x}({\bf p}_2))-n_m^2}\\
\nonumber  \hspace{-1.7ex} & = & \hspace{-1.7ex} \frac{1}{q^{n_m^2}}\left(\sum_{{\bf p}_1,\,w({\bf p}_1)=\alpha'}q^{d(\underline{x}({\bf p}_1))}\right)
\left(\sum_{{\bf p}_2,\,w({\bf p}_2)=\alpha''}q^{d(\underline{x}({\bf p}_2))}\right).
\end{eqnarray}
Thus, from \eqref{KarMin_jed} we obtain
we obtain
\begin{equation} \label{KarProdA_jed}
\chi_W^\alpha(q)=\frac{(q)_{n_m}}{q^{n_m^2}} \chi_{\g',W(\Lambda_0)}^{\alpha'}(q)
\chi_{\g'',W(\Lambda_0)}^{\alpha''}(q),
\end{equation}
where $\chi_{\g',W(\Lambda_0)}^{\alpha'}(q)$ and
$\chi_{\g'',W(\Lambda_0)}^{\alpha''}(q)$ are character formulas for Feigin-Stoyanovsky's type subspaces $W(\Lambda_0)$ for $\g'$ and $\g''$, respectively.
Formulas \eqref{kar1r_jed} and \eqref{kar2r_jed}
give
\begin{equation} \label{kar02_jed}
\chi_{W(\Lambda_0)}^{n_1,\dots,n_\ell}(q)=\frac{ q^{\sum_{i=1}^\ell n_i^2-\sum_{i=1}^{\ell-1} n_i n_{i+1}}  (q)_{n_m}  }
{(q)_{n_1}(q)_{n_2-n_1}\cdots(q)_{n_m-n_{m-1}}
(q)_{n_m-n_{m+1}} \cdots (q)_{n_{\ell-1}-n_\ell} (q)_{n_\ell}}.\end{equation}

In other cases, when $1\leq k \leq \ell$, the reasoning is similar, one only needs
to slightly modify definitions of $\underline{x}({\bf p})$, $\underline{x}({\bf p}_1)$ and $\underline{x}({\bf p}_2)$ by setting
$$r_1= r_1'=1+\theta(k-i_1),\, r_1''=1,$$
if $1\leq k \leq m$, or
$$r_1= r_1''=1+\theta(j_1-k),\, r_1'=1,$$
if $m< k \leq \ell$. In the first case, $\underline{x}({\bf p})$, $\underline{x}({\bf p}_1)$ and $\underline{x}({\bf p}_2)$ are
the smallest monomials of paths $\bf p$, ${\bf p}_1$,${\bf p}_2$, that satisfy difference and initial conditions for $L(\Lambda_k)$, $L(\Lambda_k)$
and $L(\Lambda_0)$, respectively. In the second case, these are
the smallest monomials of paths $\bf p$, ${\bf p}_1$, ${\bf p}_2$ that satisfy difference and initial conditions for $L(\Lambda_k)$, $L(\Lambda_0)$
and $L(\Lambda_k)$, respectively.
Like in $\eqref{KarProdA_jed}$, for $1\leq k \leq m$
we have $$\chi_{W(\Lambda_k)}^\alpha(q)=\frac{(q)_{n_m}}{q^{n_m^2}} \chi_{\g',W(\Lambda_k)}^{\alpha'}(q)
\chi_{\g'',W(\Lambda_0)}^{\alpha''}(q),$$
while for $m< k \leq \ell$ we have $$\chi_{W(\Lambda_k)}^\alpha(q)=\frac{(q)_{n_m}}{q^{n_m^2}} \chi_{\g',W(\Lambda_0)}^{\alpha'}(q)
\chi_{\g'',W(\Lambda_k)}^{\alpha''}(q).$$

\begin{tm} \label{KarAm_tm}
For $0\leq n_1\leq\dots\leq n_m\geq\dots\geq n_\ell\geq 0$,
\begin{equation}
\label{KarAm_jed}
\chi_{W(\Lambda_k)}^{n_1,\dots,n_\ell}(q)=\frac{q^{\sum_{i=1}^\ell n_i^2-\sum_{i=1}^{\ell-1} n_i n_{i+1}+n_k}   (q)_{n_m}  }
{(q)_{n_1}(q)_{n_2-n_1}\cdots(q)_{n_m-n_{m-1}}
(q)_{n_m-n_{m+1}} \cdots (q)_{n_{\ell-1}-n_\ell} (q)_{n_\ell}}.\end{equation}
\end{tm}

\subsection{Recurrence relations}
\label{ARR_sect}

%

We say that a monomial $$\underline{x}=x_{i_n j_n}(-r_n)\cdots x_{i_1 j_1}(-r_1)\in\C[\Gamt^-]$$ satisfies IC$_{ij}$ if  either $r_1\geq 2$ or
$r_1=1$ and $i_1\geq i$, $j_1\leq j$. We say that a monomial $\underline{x}$ satisfies IC$_0$ if  $r_1\geq 2$. Denote by
\begin{eqnarray*}\label{bij_def}
\baza_{ij} \hspace{-1.7ex} & = & \hspace{-1.7ex} \{\underline{x}\in\C[\Gamt^-] \,|\, \underline{x} \ \textrm{satisfies DC and IC}_{ij}\},\\
\baza_0 \hspace{-1.7ex} & = & \hspace{-1.7ex} \{\underline{x}\in\C[\Gamt^-] \,|\, \underline{x} \ \textrm{satisfies DC and IC}_0\}.
\end{eqnarray*}
Note that
\begin{equation} \label{KarRub_jed}
\left\{\begin{array}{l}
\baza_{W(\Lambda_i)}=\baza_{i+1,\ell},\ \textrm{for}\ i=1,\dots,m-1,\\
\baza_{W(\Lambda_j)}=\baza_{1,j-1},\ \textrm{for}\  j=m+1,\dots,\ell,\\ \baza_{W(\Lambda_0)}=\baza_{1,\ell},\quad \baza_{W(\Lambda_m)}=\baza_0.
\end{array}\right.
\end{equation}

The following lemma is a direct consequence of difference and initial conditions:

\begin{lm} \label{bazrek_lm}
$(i)$ \quad
Let $\underline{x}\in\baza_{ij}$; factorize $\underline{x}=\underline{x}_2 \underline{x}_1$ so that $\underline{x}_1$ contains all
elements of degree $-1$ and  $\underline{x}_2$ contains elements
of lower degree. Let $\underline{x}_1=x_{i_n j_n}(-1)\cdots x_{i_1 j_1}(-1)$.
Then $i\leq i_1\leq \dots \leq i_n \leq m \leq j_n \leq \dots \leq j_1 \leq j$.\\
$(ii)$ \quad $\underline{x}\in\baza_0$ if and only if $\underline{x}^+\in\baza_{1   \ell}$.
\end{lm}

For $\alpha\in Q$, define $\baza_{ij}^\alpha$,
$\baza_{0}^\alpha$, $\chi_{ij}^\alpha(q)$  and $\chi_{0}^\alpha(q)$ like
we did before.

\begin{prop} Let $\alpha=n_ 1 \alpha_1+\dots+n_\ell \alpha_\ell$, where $0\leq n_1\leq\dots\leq n_m\geq\dots\geq n_\ell\geq 0$. Then
\begin{eqnarray}
\label{k0r_jed}
 \chi_{0}^\alpha(q) \hspace{-1.7ex} & = & \hspace{-1.7ex} q^{n_m} \chi_{1\ell}^\alpha(q),\\
\label{kijr_jed}
\chi_{ij}^\alpha(q) \hspace{-1.7ex} & = & \hspace{-1.7ex} \chi_{i+1,j}^\alpha(q) + \chi_{i,j-1}^\alpha(q) - \chi_{i+1,j-1}^\alpha(q) + \\
\nonumber & &  \hspace{-1.7ex} q \chi_{i+1,j-1}^{\alpha-(ij)}(q) - q^{n_m}
\chi_{1,\ell}^{\alpha-(ij)}(q) + \\
\nonumber & &  \hspace{-1.7ex} q^{n_m} \left(\chi_{1,j-1}^{\alpha-(ij)}(q) -  \chi_{i+1,\ell}^{\alpha-(ij)}(q) -\chi_{i+1,j-1}^{\alpha-(ij)}(q)\right),\qquad
\textrm{for}\ i,j\neq m,\\
\chi_{im}^\alpha(q) \hspace{-1.7ex} & = & \hspace{-1.7ex} \chi_{i+1,m}^\alpha(q) + q \chi_{i+1,m}^{\alpha-(im)}(q),\\
\chi_{mj}^\alpha(q) \hspace{-1.7ex} & = & \hspace{-1.7ex} \chi_{m,j-1}^\alpha(q)
+ q \chi_{m,j-1}^{\alpha-(mj)}(q),\\
\label{kmr_jed}
\chi_{mm}^\alpha(q) \hspace{-1.7ex} & = & \hspace{-1.7ex} \chi_{0}^\alpha(q)
+ q^{n_m} \chi_{0}^{\alpha-(mm)}(q).
\end{eqnarray}
\end{prop}

\begin{dokaz}
To prove the first relation note that if $w(\underline{x})=\alpha$ then $d(\underline{x}^-)=n_m+d(\underline{x})$. The relation now follows from
Lemma \ref{bazrek_lm}.

We also prove the second relation; the others are proved in a similar manner.
Let $\underline{x}=x_{i_{n_m} j_{n_m}}(-r_{n_m})\cdots x_{i_1 j_1}(-r_1)\in\C[\Gamt^-]$.
If $r_1\geq 2$ or $r_1=1$ and $(i_1,j_1)\neq (i,j)$ then $\underline{x}\in\baza_{ij}^\alpha$ if and only if
$\underline{x}\in \baza_{i+1,j}^\alpha \cup \baza_{i,j-1}^\alpha$. Note also that $\baza_{i+1,j}^\alpha \cap \baza_{i,j-1}^\alpha=\baza_{i+1,j-1}^\alpha$.
This gives the first row on the right hand side of \eqref{kijr_jed}.\\
Assume $r_1=1$ and $(i_1,j_1)=(i,j)$. Set $\underline{x}_2=x_{i_{n_m} j_{n_m}}(-r_{n_m})\cdots x_{i_2 j_2}(-r_2)$.
If $r_2=1$, then, by Lemma \ref{bazrek_lm}, $\underline{x} \in \baza_{i,j}^{\alpha}$ if and only if
$\underline{x}_2\in \baza_{i+1,j-1}^{\alpha-(ij)}\setminus \baza_{0}^{\alpha-(ij)}$. Together with \eqref{k0r_jed},
this gives the second row on the right hand side of \eqref{kijr_jed}.\\
If $r_2\geq 2$, then, by difference conditions, $\underline{x} \in \baza_{i,j}^{\alpha}$
if and only if $r_2\geq 3$ or $r_2=2$ and $i_2>i$ or $j_2<j$. This is equivalent to saying that
$\underline{x}_2^+ \in \baza_{1,j-1}^{\alpha-(ij)} \cup \baza_{i+1,\ell}^{\alpha-(ij)}$.
Note also that $\baza_{1,j-1}^{\alpha-(ij)} \cap \baza_{i+1,\ell}^{\alpha-(ij)} = \baza_{i+1,j-1}^{\alpha-(ij)}$. This gives
the last row on the right hand side of \eqref{kijr_jed}.
\end{dokaz}

\begin{tm} For $\omega=\omega_1$ or $\omega=\omega_\ell$, the solution of the
system of recursions  \eqref{k0r_jed}--\eqref{kmr_jed} is given by formulas
\eqref{kar1r_jed} and \eqref{kar2r_jed}.
For $\omega=\omega_m$, $1<m<\ell$, the solution
of \eqref{k0r_jed}--\eqref{kmr_jed} is given by
\begin{equation} \label{RekRjA_jed}
\chi_{ij}^\alpha(q)=\frac{q^{\sum_{t=1}^\ell n_t^2-\sum_{t=1}^{\ell-1} n_t n_{t+1}}
\left(q^{n_{i-1}+n_{j+1}}+q^{n_m}\frac{(1-q^{n_{i-1}})(1-q^{n_{j+1}})}{1-q^{n_m}}  \right)   (q)_{n_m}   } {(q)_{n_1}(q)_{n_2-n_1}\cdots(q)_{n_m-n_{m-1}}
(q)_{n_m-n_{m+1}} \cdots (q)_{n_{\ell-1}-n_\ell} (q)_{n_\ell}},
\end{equation}
where we set $n_0=n_{\ell+1}=0$.
\end{tm}

\begin{dokaz}
For $\omega=\omega_1$ or $\omega=\omega_\ell$, the claim follows from
\eqref{KarRub_jed}.
Let $\omega=\omega_m$, $1<m<\ell$. If $i=1$ or $j=\ell$, formula \eqref{RekRjA_jed} is exactly the character formula
for the corresponding Feigin-Stoyanovsky's type subspace (see \eqref{KarRub_jed}).

If $i>1$ and $j<\ell$, we find $\chi_{ij}^\alpha(q)$ similarly to the way
we have computed characters in the previous subsection. We use the same notation as in subsection \ref{Al2_sect}.
For a path $\bf p$, define monomials $\underline{x}({\bf p})$, $\underline{x}({\bf p}_1)$ and $\underline{x}({\bf p}_2)$ by \eqref{MinDegRek_jed},
but this time, instead of \eqref{MinMonPocDeg_jed}, we set
\begin{eqnarray*}
 & r_1=1+\max\{\theta(i-1-i_1),\theta(j_1-j-1)\}, & \\
& r_1'=1+\theta(i-1-i_1),\quad  r_1''=1+\theta(j_1-j-1). &
\end{eqnarray*}
The difference from the previous case is that now formulas \eqref{MinDeg_jed} and \eqref{DegProdA_jed} fail for a
path $\bf p$ that starts with a color $(i_1 j_1)$ such that $i_1< i$ and $j_1>j$.
For such path, we have $r_1=r_1'=r_1''=2$, so \eqref{MinDeg_jed} does not hold for $t=1$.
This means that formula \eqref{DegProdA_jed} calculates $d(\underline{x}({\bf p}))$ as if
$r_1=3$ instead of $r_1=2$, and the difference between the calculated and the actual degree for monomials of such path is equal to $l({\bf p})=n_m$.

Although we cannot use \eqref{KarProdA_jed} to calculate $\chi_{ij}^\alpha(q)$,
we can ``repair'' the wrong character formula obtained from \eqref{KarProdA_jed}
by recalculating degrees of monomials that start with a color $(i_1 j_1)$ such that $i_1< i$, $j_1>j$. Let
\begin{eqnarray*}
\mathfrak{C}_{ij}^\alpha \hspace{-1.7ex} & = & \hspace{-1.7ex} \{x_{i_{n_m} j_{n_m}}(-r_{n_m})\cdots
x_{i_1 j_1}(-r_1) \in \baza_{ij}^\alpha \st i_1< i, j_1>j\},\\
\mathfrak{D}_{ij}^\alpha \hspace{-1.7ex} & = & \hspace{-1.7ex} \{x_{i_{n_m} j_{n_m}}(-r_{n_m})\cdots
x_{i_1 j_1}(-r_1) \in \baza_{1\ell}^\alpha \st r_1=1,i_1< i, j_1>j\}.
\end{eqnarray*}
Denote by $\chi_{\mathfrak{C}_{ij}}^\alpha(q)$ and $\chi_{\mathfrak{D}_{ij}}^\alpha(q)$ the corresponding graded cardinalities.
By observations above, we have
\begin{equation}
\label{KarRepair_jed} \nonumber
\chi_{ij}^\alpha(q) =\frac{(q)_{n_m}}{q^{n_m^2}} \chi_{\g',W(\Lambda_{i-1})}^{\alpha'}(q)
\chi_{\g'',W(\Lambda_{j+1})}^{\alpha''}(q) -
\left(1-q^{n_m}\right)\chi_{\mathfrak{C}_{ij}}^\alpha(q)
\end{equation}
Since $\underline{x}\in\mathfrak{C}^\alpha_{ij}$ if and only if  $\underline{x}^{+r}\in\mathfrak{D}^\alpha_{ij}$, for some $r\in\N$, we have
\begin{equation}\label{KriviRek1_jed} \nonumber
\chi_{\mathfrak{C}_{ij}}^\alpha(q)=\frac{q^{n_m}}{1-q^{n_m}}\chi_{\mathfrak{D}_{ij}}^\alpha(q).
\end{equation}
Furthermore, since $\mathfrak{D}_{ij}^\alpha=\baza_{1\ell}^\alpha \setminus \left( \baza_{1j}^\alpha \cup \baza_{i\ell}^\alpha \right)$ and
$\baza_{1j}^\alpha \cap \baza_{i\ell}^\alpha=\baza_{ij}^\alpha$, we have
\begin{equation}\label{KriviRek2_jed} \nonumber
\chi_{\mathfrak{D}_{ij}}^\alpha(q)=\chi_{1\ell}^\alpha(q) -
\chi_{i\ell}^\alpha(q)-\chi_{1j}^\alpha(q)
+\chi_{ij}^\alpha(q).
\end{equation}
Consequently
\begin{eqnarray*}
\chi_{ij}^\alpha(q) \hspace{-1.7ex} & = & \hspace{-1.7ex} \frac{(q)_{n_m}}{q^{n_m^2}} \chi_{\g',W(\Lambda_{i-1})}^{\alpha'}(q)
\chi_{\g'',W(\Lambda_{j+1})}^{\alpha''}(q) -  \\
 & & \hspace{-1.7ex}
q^{n_m}
\left(\chi_{1\ell}^\alpha(q) -
\chi_{i\ell}^\alpha(q)-\chi_{1j}^\alpha(q)+\chi_{ij}^\alpha(q)\right).
\end{eqnarray*}
Formula \eqref{RekRjA_jed} now follows from \eqref{kar1r_jed} and \eqref{kar2r_jed}, and Theorem \ref{KarAm_tm}.
\end{dokaz}

\section{The case $D_4$}

\label{d4_sect}

\subsection{Character formula for $W(\Lambda_0)$}

Let ${\mathfrak g}$ be a simple finite-dimensional Lie algebra of type $D_\ell$. The minuscule fundamental weights are $\omega_1, \omega_{\ell-1}, \omega_\ell$.
Fix a minuscule weight $\omega=\omega_1$.
The set of colors is $\Gamma=\{\underline{2},\dots,\underline{\ell},\ell,\dots,2\}$, where
\begin{equation}\label{gam_jed}
\left\{\begin{array}{rcl}
\underline{2} & = & \alpha_1,  \\
 & \vdots & \\
\underline{\ell} & = & \alpha_1+\dots+\alpha_{\ell-1}, \\
\ell & = & \alpha_1+\dots+\alpha_{\ell-2} + \alpha_{\ell}, \\
\ell-1 & = &  \alpha_1+\dots+\alpha_{\ell-2} +\alpha_{\ell-1}+ \alpha_{\ell}, \\
\ell-2 & = &  \alpha_1+\dots+\alpha_{\ell-3}+2\alpha_{\ell-2} +\alpha_{\ell-1}+ \alpha_{\ell}, \\
\ell-3 & = &  \alpha_1+\dots+\alpha_{\ell-4}+2\alpha_{\ell-3}+2\alpha_{\ell-2} +\alpha_{\ell-1}+ \alpha_{\ell}, \\
 & \vdots & \\
2 & = & \alpha_1+2\alpha_2+\dots++2\alpha_{\ell-2} +\alpha_{\ell-1}+ \alpha_{\ell}.
\end{array}\right.
\end{equation}
Define an order on $\Gamma$ by setting:
$2>\dots>\ell>\underline{\ell}>\dots>\underline{2}$.
Like in the previous section, this induces the order on $\Gamt$, and we assume that monomials from $\C[\Gamt]$ are sorted ascendingly from left to right.

Let $L(\Lambda_k)$, $k=0,1,\ell-1$ or $\ell$, be a standard $\gt$-module of level 1, and set $W=W(\Lambda_k)$.
Define an energy function $E:\Gamma\times\Gamma\to\{0,1,2\}$ by
\begin{equation}\label{ener2_jed}
E(\gamma',\gamma)=\left\{\begin{array}{ll}
0,  &  (\gamma',\gamma)=(\underline{2},2),\\
1,  &  \gamma'<\gamma,(\gamma',\gamma)\neq(\underline{2},2)\ \textrm{or}\ (\gamma',\gamma)=(\ell,\underline{\ell}),\\
2,  &  \gamma'\geq\gamma,(\gamma',\gamma)\neq(\ell,\underline{\ell}).
\end{array}\right.
\end{equation}
We say that a monomial $\underline{x}=x_{\gamma_n}(-r_n)\cdots x_{\gamma_1}(-r_1)\in\C[\Gamt^-]$ satisfies \emph{difference conditions}, or DC for short, if
$r_{t+1}-r_t \geq E(\gamma_{t+1},\gamma_t)$.
We say that $\underline{x}$ satisfies
{\em initial conditions} for $L(\Lambda_k)$, or IC for short, if either $r_1\geq 2$ or
$r_1=1$ and $\gamma_1\in\{\underline{2},\dots,\underline{\ell-1},\ell\}$, for $k=\ell-1$, or $\gamma_1\in\{\underline{2},\dots,\underline{\ell}\}$, for $k=\ell$, or
$\gamma_1\in\{\underline{2},\dots,\underline{\ell},\ell,\dots,2\}$, for $k=0$.
As before, define the set $\baza_W$ by \eqref{baza_def}.

\begin{tm}[{[B],[P2]}]
The set
$\{\underline{x} v_{\Lambda_k} \,|\, \underline{x}\in\baza_W \}$
is a basis of $W$.
\end{tm}

From now on we assume that the algebra $\g$ is of  type $D_4$; $\Gamma=\{\underline{2},\underline{3},\underline{4},4,3,2\}$.
Like in the previous section, we define weight and degree of monomials, and the character $\chi_W(z_1,z_2,z_3,z_4,q)$ of $W$. Furthermore,
for $n_1,n_2,n_3,n_4\geq 0$  set $\alpha=n_1 \alpha_1+ n_2 \alpha_2 + n_3 \alpha_3 + n_4 \alpha_4$ and define sets $\baza_W^\alpha$ and
formal series $\chi_W^\alpha (q)$ as before.

Obviously, $\chi_W^\alpha(q)=0$ unless $\alpha$ can be written in the form
\begin{equation}\label{alpha_jed}
\alpha=m_{\underline{2}} \underline{2} + m_{\underline{3}} \underline{3} + m_{\underline{4}} \underline{4} +
m_4 4 + m_3 3 + m_2 2,
\end{equation} for some $m_{\underline{2}}, m_{\underline{3}}, m_{\underline{4}}, m_4, m_3, m_2\in\Z_{\geq 0}$. Set
\begin{equation}\label{nula_jed}
\underline{0}=2+\underline{2}=3+\underline{3}=4+\underline{4}=2\alpha_1+2\alpha_2+\alpha_3+\alpha_4;
\end{equation}
then \eqref{alpha_jed} is equivalent to
\begin{equation} \label{alpha2_jed}
\alpha=m_{\underline{2}} \underline{2} + m_{\underline{3}} \underline{3} + m_{\underline{4}} \underline{4} +
m_{\underline{0}} \underline{0},
\end{equation}
 where $m_{\underline{2}}, m_{\underline{3}}, m_{\underline{4}}, m_{\underline{0}}\in\Z$, such that
\begin{equation} \label{AlfaUvijet_jed}
m_{\underline{0}}\geq -\theta(-m_{\underline{2}}) m_{\underline{2}} - \theta (-m_{\underline{3}}) m_{\underline{3}} - \theta(-m_{\underline{4}}) m_{\underline{4}}.
\end{equation}
By \eqref{gam_jed} and \eqref{nula_jed} we have
\begin{equation}\label{m2n_jed}
m_{\underline{2}}  =   n_1-n_2,\
m_{\underline{3}}  =   n_2-n_3-n_4,\
m_{\underline{4}}  =   n_3-n_4,\
m_{\underline{0}}  =   n_4.
\end{equation}
Therefore condition \eqref{AlfaUvijet_jed} is equivalent to
\begin{equation} \label{AlfaUvijet2_jed}
\left\{\begin{array}{ll}
n_1-n_2+n_4\geq 0, & n_1-n_2+n_3\geq 0,\\
n_2-n_3\geq 0, & n_1-n_3\geq 0,\\
n_3\geq 0, & n_2-n_4\geq 0,\\
n_4\geq 0, & n_1-n_4\geq 0.\\
\end{array}\right.
\end{equation}

We first consider the case $W=W(\Lambda_0)$. The other cases will be considered in the next subsection.

Set $\Gamma'=\{\underline{2},\underline{4},4,2\}$, $\Gamma''=\{\underline{3},3\}$. Define
\begin{eqnarray*}
\baza_{\Gamma'} \hspace{-1.7ex} & = & \hspace{-1.7ex} \{x_{\gamma_n}(-r_n)\cdots x_{\gamma_1}(-r_1) \in \baza_{W(\Lambda_0)} \st
\gamma_i\in\Gamma',\, i=1,\dots,n \}, \\
\baza_{\Gamma''} \hspace{-1.7ex} & = & \hspace{-1.7ex} \{x_{\gamma_n}(-r_n)\cdots x_{\gamma_1}(-r_1) \in \baza_{W(\Lambda_0)} \st
\gamma_i\in\Gamma'',\, i=1,\dots,n \},\end{eqnarray*}
$\baza_{\Gamma'}^\alpha=\baza_{\Gamma'}\cap \baza_{W(\Lambda_0)}^\alpha$ and $\baza_{\Gamma''}^\alpha=\baza_{\Gamma''}\cap \baza_{W(\Lambda_0)}^\alpha$.
Define formal series $\chi_{\Gamma'}^\alpha$ and $\chi_{\Gamma''}^\alpha$ in the obvious way.

By setting
\begin{equation} \label{D4ident1_jed}
\underline{2}=(22),\ \underline{4}=(23),\ 4=(12),\ 2=(13),
\end{equation}
 we identify the set $\Gamma'$ with the set of colors from the case
$A_3, \omega=\omega_2$ (see Section \ref{Al2_sect}).
Since the energy functions agree with this identification, and since in both cases we have the same relations between colors:
$$\underline{2}+2=\underline{4}+4,\quad\textrm{i.e.}\quad (22)+(13)=(23)+(12),$$
we conclude that the sets of monomials satisfying difference and initial conditions coincide. Therefore we can deduce a formula
for $\chi_1^{n_1,n_2,n_3,n_4}(q)$ from the character formula for $W(\Lambda_0)$ for $A_3, \omega=\omega_2$. Let
\begin{equation}
\label{alph1_jed}
\alpha'=n_1 \alpha_1+ n_2 \alpha_2 + n_3 \alpha_3 + n_4 \alpha_4=
m_{\underline{2}} \underline{2} + m_{\underline{4}} \underline{4} +
m_4 4 + m_2 2,
\end{equation}
for some $m_{\underline{2}}, m_{\underline{4}}, m_4, m_2\geq 0$.
By \eqref{gam_jed}, we have
\begin{equation} \label{n2m1_jed}
n_1=m_{\underline{2}}+m_{\underline{4}}+m_4+m_2,\
n_2=m_{\underline{4}}+m_4+2 m_2,\
n_3=m_{\underline{4}}+m_2,\
n_4=m_4+m_2.
\end{equation}
Note from \eqref{retci_jed}, \eqref{stupci_jed} and
\eqref{D4ident1_jed} that parameters $n_1, n_2-n_1,n_2-n_3,n_3$ from
the case $A_3, \omega=\omega_2$, correspond to
$m_2+m_4,m_{\underline{2}}+m_{\underline{4}},
m_4+m_{\underline{2}},m_{\underline{4}}+m_2$ from the
$\Gamma'$-case, respectively. From \eqref{n2m1_jed} we see that in
the $\Gamma'$-case these parameters are equal to $n_4,
n_1-n_4,n_1-n_3,n_3$, respectively, and are independent of the
particular choice of $m_{\underline{2}}, m_{\underline{4}}, m_4,
m_2$. Moreover, $n_2=n_1+(n_2-n_1)$ from the case $A_3,
\omega=\omega_2$, corresponds to
$n_1=m_{\underline{2}}+m_{\underline{4}}+m_4+m_2$ from the
$\Gamma'$-case. Hence, from  character formula \eqref{kar02_jed} for
$W(\Lambda_0)$ for the case $A_3, \omega=\omega_2$, we get
\begin{equation}\label{kar1_jed}
\chi_{\Gamma'}^{\alpha'} (q)=\frac{q^{n_4^2+n_1^2+n_3^2-n_4 n_1-n_1 n_3}}{(q)_{n_3}(q)_{n_1-n_3}(q)_{n_1-n_4}
(q)_{n_4}}.
\end{equation}

Similarly, we identify the set $\Gamma''$ with the set of
colors from the case $A_2, \omega=\omega_2$ (see Section \ref{Al1_sect}):
\begin{equation} \label{D4ident2_jed}
\underline{3}=(2),3=(1).
\end{equation}
Let
\begin{equation}
\label{alph2_jed}
\alpha''=n_1 \alpha_1+ n_2 \alpha_2 + n_3 \alpha_3 + n_4 \alpha_4=
m_{\underline{3}} \underline{3} + m_3 3,
\end{equation}
for some $m_{\underline{3}}, m_3\geq 0$.
Then, by \eqref{gam_jed},
\begin{equation} \label{n2m2_jed}
n_1=m_{\underline{3}}+m_3,\
n_2=m_{\underline{3}}+m_3,\
n_3=m_3,\
n_4=m_3.
\end{equation}
From \eqref{stupci_jed} and \eqref{D4ident2_jed} we see that the parameters $n_1, n_2-n_1, n_2$ from the case $A_2, \omega=\omega_2$
correspond to $n_4, n_1-n_4, n_1$ from the $\Gamma''$-case. Hence, from character formula \eqref{kar1r_jed} for $W(\Lambda_0)$ for $A_2, \omega=\omega_2$, we get
\begin{equation}\label{kar2_jed}
\chi_{\Gamma''}^{\alpha''}(q)=\frac{q^{n_1^2+n_4^2-n_1 n_4}}{(q)_{n_4}(q)_{n_1-n_4}}.
\end{equation}

The following procedure gives us a way to obtain a character formula for $W$
from formulas \eqref{kar1_jed} and \eqref{kar2_jed}.
Set $\Gamma^e=\Gamma\cup\{\tilde{2},\tilde{4}\}$ and $\Gamt^e=\Gamt\cup \{x_{\gamma}(-r)\st \gamma\in \{\tilde{2},\tilde{4}\},\,r\in\Z\}$.
Define $\tilde{2}>2>3>4>\tilde{4}>\underline{4}>\underline{3}>\underline{2}$, and define the order on $\Gamt^e$ accordingly.
Let $\underline{x}_1\in\baza_{\Gamma'}$, $\underline{x}_2\in\baza_{\Gamma''}$.
Denote by $\underline{x}_3\in\C[\Gamt^e]$ a monomial obtained from $\underline{x}_1$ by replacing every pair $x_{\underline{2}}(-r) x_2(-r)$
 with a pair $x_{\tilde{2}}(-r-1)x_{\tilde{2}}(-r)$, and every pair $x_{4}(-r-1) x_{\underline{4}}(-r)$ with a pair $x_{\tilde{4}}(-r-1)x_{\tilde{4}}(-r)$.
Set $\underline{y}=\underline{x}_3^{+\triangledown} \underline{x}_2^{+\triangledown}$,  and reorder variables so that they are
sorted ascendingly from  left to right. Set $\underline{z}=\underline{y}^{-\triangledown}$; note that pairs $x_{\gamma}(-r-1)x_{\gamma}(-r)$, $\gamma \in\{
\tilde{2},\tilde{4}\}$, from $\underline{x}_3$ correspond to pairs $x_{\gamma}(-r'-1)x_{\gamma}(-r')$ from $\underline{z}$.
Let $\underline{x}\in \C[\Gamt]$ be a monomial obtained from $\underline{z}$ by replacing every pair $x_{\tilde{2}}(-r-1)x_{\tilde{2}}(-r)$
inside $\underline{z}$ with a pair $x_{\underline{2}}(-r)x_2(-r)$, and every pair $x_{\tilde{4}}(-r-1)x_{\tilde{4}}(-r)$
with a pair $x_4(-r-1)x_{\underline{4}}(-r)$.


\begin{prop} \label{KarKomb_prop}
Let $\underline{x},\underline{x}_1,\underline{x}_2$ be as above. Then
$\underline{x}$ satisfies difference and initial conditions. Conversely,
every monomial that satisfies difference and initial conditions can be obtained
in this way.
\end{prop}

\begin{dokaz}
Let $\underline{x}=x_{\gamma_n}(-r_n)\cdots x_{\gamma_1}(-r_1)$. For
$t=1,\dots,n-1$, consider factors $x_{\gamma_t}(-r_t)$ and
$x_{\gamma_{t+1}}(-r_{t+1})$. If $\gamma_{t},\gamma_{t+1}\in\Gamma'$
or $\gamma_{t},\gamma_{t+1}\in\Gamma''$ then these two factors
obviously satisfy difference conditions since they come from the two
neighbouring factors inside $\underline{x}_1$ or $\underline{x}_2$,
respectively, and the above procedure did not change the difference
between their degrees.

If $\gamma_t\in\Gamma'$ and  $\gamma_{t+1}\in\Gamma''$ or $\gamma_t\in\Gamma''$ and  $\gamma_{t+1}\in\Gamma'$, then either $r_{t+1}-r_t=1$
and $\gamma_{t+1}<\gamma_t$, or $r_{t+1}-r_t\geq 2$, which means that difference conditions are again satisfied.

Conversely, let $\underline{x}=x_{\gamma_n}(-r_n)\cdots x_{\gamma_1}(-r_1)\in\C[\Gamt^-]$ be a monomial that satisfies difference and initial conditions.
Let $\underline{x}_1$ and $\underline{x}_2$
be monomials obtained by the reverse procedure.
The claim will follow from the following simple observations that can be proved inductively from \eqref{ener2_jed}:\\
\noindent $(i)$\quad If $\gamma_{t},\gamma_{t+1},\dots,\gamma_{t+s}\in\Gamma'$, then
$r_{t+s}-r_t\geq s-1$. Moreover, $r_{t+s}-r_t=s-1$ if and only if $s$ is odd and $(\gamma_t,\dots,\gamma_{t+s})=(2,\underline{2},\dots,2,\underline{2}).$\\
\noindent $(ii)$\quad If $\gamma_{t},\gamma_{t+1},\dots,\gamma_{t+s}\in\Gamma''$, then
$r_{t+s}-r_t\geq s+q$, where $q=\#\{0\leq i\leq s-1 \st  \gamma_{t+i}<\gamma_{t+i+1}\}$.

First we show that $\underline{x}_2$ satisfies difference conditions. Let $x_\gamma(-r), x_{\gamma'}(-r')$ be two neighbouring factors inside $\underline{x}_2$.
Assume that $x_{\gamma_t}(-r_t)$ and $x_{\gamma_{t+s}} (-r_{t+s})$  are the corresponding factors inside $\underline{x}$. If $s=1$, then it is obvious that
$x_\gamma(-r)$ and $x_{\gamma'}(-r')$ satisfy difference condition.
If $s>1$, then $\gamma_{t+1},\dots,\gamma_{t+s-1}\in\Gamma'$ and $r'-r=r_{t+s}-r_{t}-s+1$.
We need to show that either $r_{t+s}-r_{t}\geq s+1$ or $r_{t+s}-r_{t}= s$ and $\gamma>\gamma'$. By $(i)$ and \eqref{ener2_jed} we have:
\begin{equation} \label{rac4_jed}
r_{t+1}-r_t\geq 1,\ r_{t+s-1}-r_{t+1}\geq s-3,\ r_{t+s}-r_{t+s-1}\geq 1.
\end{equation}
Moreover, if $r_{t+s-1}-r_{t+1}=s-3$, by $(i)$ we have $\gamma_{t+s-1}=\underline{2}$ and therefore $r_{t+s}-r_{t+s-1}\geq 2$. Hence,
$r_{t+s}-r_{t}\geq s$. Assume that $r_{t+s}-r_{t}= s$ (this is the case when $r'=r+1$). Then \eqref{rac4_jed} and \eqref{ener2_jed} imply
\begin{equation} \label{rac3_jed}
\gamma_t>\gamma_{t+1}.
\end{equation} Hence
\begin{equation} \label{rac2_jed}
\gamma_{t+1}\neq 2
\end{equation}
 and $r_{t+s-1}-r_{t+1}=s-2$, $r_{t+s}-r_{t+s-1}= 1$. By \eqref{ener2_jed}, we see that
\begin{equation} \label{rac1_jed}
\gamma_{t+s-1}>\gamma_{t+s}.
\end{equation} If $\gamma_{t+1}=\underline{2}$ then $r_{t+s-1}-r_{t+2}=s-4$, so from $(i)$ we see that $\gamma_{t+s-1}=\underline{2}$.
But this is in contradiction with \eqref{rac1_jed}. If $\gamma_{t+s-1}=2$, then $r_{t+s-2}-r_{t+1}=s-4$. By $(i)$, this implies $\gamma_{t+1}=2$
which is in contradiction with \eqref{rac2_jed}. So, if $r_{t+s}-r_{t}= s$, then $\gamma_{t+1},\gamma_{t+s-1}\in\{\underline{4},4\}$.
By \eqref{rac3_jed} and \eqref{rac1_jed} we conclude that $\gamma=\gamma_t>\gamma_{t+s}=\gamma'$. Therefore  $x_\gamma(-r)$ and $x_{\gamma'}(-r-1)$
satisfy difference conditions.

In the same way we show that $\underline{x}_1$ satisfies difference conditions. Let $x_\gamma(-r)$ and $x_{\gamma'}(-r')$ be two neighbouring factors
inside $\underline{x}_1$. Assume that $x_{\gamma_t}(-r_t)$ and $x_{\gamma_{t+s}} (-r_{t+s})$  are the corresponding factors inside $\underline{x}$.
Again, if $s=1$, the claim is obvious. Assume $s>1$.
Then $\gamma_{t+1},\dots,\gamma_{t+s-1}\in\Gamma''$ and $r'-r=r_{t+s}-r_{t}-s+1$.
We need to show that either $r_{t+s}-r_{t}\geq s+1$ or $r_{t+s}-r_{t}= s$ and $\gamma>\gamma'$.
By $(ii)$ and \eqref{ener2_jed} we have:
$$r_{t+1}-r_t\geq 1,\ r_{t+s-1}-r_{t+1}\geq s+q-2,\ r_{t+s}-r_{t+s-1}\geq 1,$$
where $q$ is defined in $(ii)$. Therefore
$r_{t+s}-r_{t}\geq s+q$. If $q=0$ and $r_{t+s}-r_{t} = s$, by \eqref{ener2_jed}, we must have
$\gamma_t>\gamma_{t+1}>\dots>\gamma_{t+s-1}>\gamma_{t+s}$.
Hence $\gamma=\gamma_t>\gamma_{t+s}=\gamma'$ and $r'=r+1$. We conclude that, in this case, $x_\gamma(-r)$ and $x_{\gamma'}(-r-1)$ satisfy difference conditions.

By using similar arguments we can show that $\underline{x}_1,\underline{x}_2\in\C[\Gamt^-]$, i.e. that  factors of $\underline{x}_1$ and $\underline{x}_2$
have negative degrees. Hence, $\underline{x}_1\in \baza_{\Gamma'}$ and $\underline{x}_2\in\baza_{\Gamma''}$.
\end{dokaz}

Let $w(\underline{x}_1)=n_1' \alpha_1+ n_2' \alpha_2 + n_3' \alpha_3 + n_4' \alpha_4$ and
$w(\underline{x}_2)=n_1'' \alpha_1+ n_2'' \alpha_2 + n_3'' \alpha_3 + n_4'' \alpha_4$.
From the construction we see
\begin{eqnarray}
d(\underline{x}) \hspace{-1.7ex} & = & \hspace{-1.7ex} d(\underline{x}_1) + d(\underline{x}_2) - \frac{n_1'(n_1'-1)}{2} -
\frac{n_1''(n_1''-1)}{2} + \frac{(n_1'+n_1'')(n_1'+n_1''-1)}{2} \nonumber\\
 \hspace{-1.7ex} & = & \hspace{-1.7ex} d(\underline{x}_1) + d(\underline{x}_2) + n_1' n_1''. \label{Kor_jed}
\end{eqnarray}

Fix $n_1,n_2,n_3,n_4\geq 0$ satisfying \eqref{AlfaUvijet2_jed}, and set $\alpha=n_1 \alpha_1+ n_2 \alpha_2 + n_3 \alpha_3 + n_4 \alpha_4$. Define
$m_{\underline{2}}, m_{\underline{3}}, m_{\underline{4}}, m_{\underline{0}}$ by
\eqref{m2n_jed}.
Define
\begin{equation} \label{MCrta_jed}
m'= -\theta(-m_{\underline{2}}) m_{\underline{2}} - \theta(-m_{\underline{4}}) m_{\underline{4}},\qquad
m''= - \theta (-m_{\underline{3}}) m_{\underline{3}}.
\end{equation}
Condition \eqref{AlfaUvijet2_jed} is equivalent to $m_{\underline{0}}\geq m'+m''$ (cf. \eqref{AlfaUvijet_jed}).
For $i=0,\dots,m_{\underline{0}}-m'-m''$, set
\begin{equation} \label{AlfaI_jed}
\alpha_i'  =  m_{\underline{2}} \underline{2} + m_{\underline{4}} \underline{4} + (i+m') \underline{0},\qquad
\alpha_i''  =   \alpha-\alpha_i'.
\end{equation}
By \eqref{gam_jed} and \eqref{nula_jed} we have
\begin{eqnarray*}
\alpha_i'  \hspace{-1.7ex} & = & \hspace{-1.7ex}  (n_1-n_2+n_3-n_4+2(i+m'))\alpha_1
+(n_3-n_4+2(i+m'))\alpha_2\\
& & +(n_3-n_4+i+m')\alpha_3 +(i+m')\alpha_4,\\
\alpha_i''  \hspace{-1.7ex} & = & \hspace{-1.7ex}   (n_2-n_3+n_4-2(i+m'))\alpha_1
+(n_2-n_3+n_4-2(i+m'))\alpha_2\\
& & +(n_4-i-m')\alpha_3 +(n_4-i-m')\alpha_4.
\end{eqnarray*}

Then, by Proposition \ref{KarKomb_prop}, \eqref{m2n_jed} and \eqref{Kor_jed},
\begin{equation} \label{KarKomb_jed}
\chi_W^\alpha=\sum_{i=0}^{n_4-m'-m''}
\chi_1^{\alpha_i'}(q) \chi_2^{\alpha_i''}(q) q^{(n_1-n_2+n_3-n_4+2(i+m'))(n_2-n_3+n_4-2(i+m'))   }.
\end{equation}
From \eqref{kar1_jed} and \eqref{kar2_jed} we obtain the following character formula:
\begin{tm}\label{KarD_tm}
\begin{eqnarray} \label{KarD_jed}
\chi_{W(\Lambda_0)}^\alpha \hspace{-1.7ex} & = & \hspace{-1.7ex} \sum_{i=0}^{n_4-m'-m''} q^{f(\alpha)} \frac{(q)_{n_1-n_2+n_3-n_4+2(i+m')}}{
(q)_{n_3-n_4+i+m'} (q)_{n_1-n_2+i+m'} (q)_{i+m'} } \\
  &  & \hspace{5ex} \cdot\frac{1}{(q)_{n_1-n_2+n_3-n_4+i+m'} (q)_{n_4-i-m'} (q)_{n_2-n_3-i-m'}    }, \nonumber
\end{eqnarray}
where
\begin{eqnarray} \label{KvadD_jed}
f(\alpha) \hspace{-1.7ex} & = & \hspace{-1.7ex} n_1^2+n_2^2+n_3^2+n_4^2-n_1 n_2 - n_2 n_3 - n_3 n_4\\
 & & \hspace{-1.7ex} -(i+m') (n_2 - n_3 + n_4 - i - m'). \nonumber
\end{eqnarray}
\end{tm}

\begin{napomena}{\em
For an algebra $\g$ of  type $D_\ell$,
let $\alpha=n_1 \alpha_1+\dots + n_\ell \alpha_\ell$ be a weight that can be written as a non-negative linear combination of
colors $\underline{2},\dots,\underline{\ell},\ell,\dots,2$. Set
$\underline{0}=2+\underline{2}=\dots=\ell+\underline{\ell}=2\alpha_1+\dots+2\alpha_{\ell-2}+\alpha_{\ell-1}+\alpha_\ell$.
Define
$m_{\underline{2}}, \dots, m_{\underline{\ell}}, m_{\underline{0}}$ like in
\eqref{alpha2_jed}. Like in \eqref{AlfaUvijet_jed}, we obtain the following condition on the coefficients
$m_{\underline{2}}, \dots, m_{\underline{\ell}}, m_{\underline{0}}$:
$$
m_{\underline{0}}\geq -\theta(-m_{\underline{2}}) m_{\underline{2}} - \dots - \theta(-m_{\underline{\ell}}) m_{\underline{\ell}}.
$$
Partition the set of colors into the sets
$$\Gamma^{(2)}=\{\underline{2},\underline{\ell},\ell,2\},\ \Gamma^{(3)}=\{\underline{3},3\},\ \dots\ ,\, \Gamma^{(\ell-1)}=\{\underline{\ell-1},\ell-1\},$$
and regard $\Gamma^{(2)}$ as a set of colors for the case $A_3,\omega=\omega_2$, and $\Gamma^{(3)},\dots\ ,\Gamma^{(\ell-1)}$
as sets of colors for the case $A_2,\omega=\omega_2$. Set $m^{(2)}=-\theta(-m_{\underline{2}}) m_{\underline{2}} -
\theta(-m_{\underline{\ell}}) m_{\underline{\ell}}$, $m^{(3)}=-\theta(-m_{\underline{3}}) m_{\underline{3}}$, \dots,
$m^{(\ell-1)}=-\theta(-m_{\underline{\ell-1}}) m_{\underline{\ell-1}},$
and let $m=m^{(2)}+\dots+m^{(\ell-1)}$.
We can apply the same procedure as before; we obtain the following character formula:
\begin{eqnarray*}
\chi_{W(\Lambda_0)}^\alpha \hspace{-1.7ex} & = & \hspace{-1.7ex} \sum_{\substack{i_2,\dots,i_{\ell-1}\geq 0
  \\ i_2+\dots+i_{\ell-1}=n_\ell-m}}  q^{f(\alpha)} \frac{(q)_{n_1-n_2+n_{\ell-1}-n_\ell+2(i_2+m^{(2)})}}{
(q)_{n_{\ell-1}-n_\ell+i_2+m^{(2)}} (q)_{n_1-n_2+i_2+m^{(2)}} (q)_{i_2+m^{(2)}} } \\
  &  & \hspace{5ex} \cdot \frac{1}{(q)_{n_1-n_2+n_{\ell-1}-n_\ell+i_2+m^{(2)}}}\prod_{j=3}^{\ell-2}
  \frac{1}{(q)_{i_j+m^{(j)}} (q)_{n_{j-1}-n_j+i_j+m^{(j)}}}  \\
  &  & \hspace{5ex} \cdot \frac{1}{(q)_{i_{\ell-1}+m^{(\ell-1)}} (q)_{n_{\ell-2}-n_{\ell-1}-n_\ell+i_{\ell-1}+m^{(\ell-1)}}},
\end{eqnarray*}
where
\begin{eqnarray*}
f(\alpha) \hspace{-1.7ex} & = & \hspace{-1.7ex} n_1^2+\dots+n_\ell^2-n_1 n_2 -\dots-n_{\ell-3}n_{\ell-2} - n_{\ell-2} n_{\ell-1} -
n_{\ell-2} n_{\ell}\\
 & & \hspace{-1.7ex} + \sum_{j=3}^{\ell-2} (i_j+m^{(j)})(n_{j-1}-n_j+i_j+m^{(j)})\\
 & & \hspace{-1.7ex} + (i_{\ell-1}+m^{(\ell-1)})(n_{\ell-2}-n_{\ell-1}-n_\ell+i_{\ell-1}+m^{(\ell-1)}).
\end{eqnarray*}
}\end{napomena}

\subsection{Character formulas for other level $1$ standard modules and recurrence relations}

For $\gamma\in\Gamma$, we say that a monomial $$\underline{x}=x_{\gamma_n}(-r_n)\cdots x_{\gamma_1}(-r_1)\in\C[\Gamt^-]$$
satisfies IC$_{\gamma}$ if  either $r_1\geq 2$ or
$r_1=1$ and either $\gamma_1\leq\gamma$ if $\gamma\neq 4$,
or $\gamma_1\in\{\underline{2},\underline{3},4\}$ if $\gamma=4$. We say that a monomial $\underline{x}$ satisfies IC$_0$ if  $r_1\geq 2$. Denote by
\begin{eqnarray*}\label{bgamma_def}
\baza_{\gamma} \hspace{-1.7ex} & = & \hspace{-1.7ex} \{\underline{x}\in\C[\Gamt^-] \,|\, \underline{x} \ \textrm{satisfies DC and IC}_{\gamma}\},\\
\baza_0 \hspace{-1.7ex} & = & \hspace{-1.7ex} \{\underline{x}\in\C[\Gamt^-] \,|\, \underline{x} \ \textrm{satisfies DC and IC}_0\}.
\end{eqnarray*}
Note that
\begin{equation} \label{KarRubD_jed}
\baza_{W(\Lambda_0)}=\baza_{2},\quad \baza_{W(\Lambda_1)}=\baza_0,\quad
\baza_{W(\Lambda_3)}=\baza_{4},\quad \baza_{W(\Lambda_4)}=\baza_{\underline{4}}.
\end{equation}

For $n_1,n_2,n_3,n_4$ satisfying
\eqref{AlfaUvijet2_jed}, set $\alpha=n_1 \alpha_1+n_2 \alpha_2+n_3 \alpha_3+n_4 \alpha_4$ and define $\baza_{\gamma}^\alpha$,
$\baza_{0}^\alpha$, $\chi_{\gamma}^\alpha(q)$  and $\chi_{0}^\alpha(q)$ as
before.

\begin{prop}  Characters $\chi_{\gamma}^\alpha(q)$, $\gamma\in\Gamma\cup\{0\}$, satisfy
the following recurrence relations:
\begin{eqnarray*}
\chi_{0}^\alpha(q) \hspace{-1.7ex} & = & \hspace{-1.7ex} q^{n_1} \chi_{2}^\alpha(q),\\
\chi_{2}^\alpha(q) \hspace{-1.7ex} & = & \hspace{-1.7ex} \chi_{3}^\alpha(q) + q^{n_1} \chi_{3}^{\alpha-2}(q) + q^{2n_1-2}
             \chi_{2}^{\alpha-2-\underline{2}}(q),\\
\chi_{3}^\alpha(q) \hspace{-1.7ex} & = & \hspace{-1.7ex} \chi_{4}^\alpha(q)+\chi_{\underline{4}}^\alpha(q)-\chi_{\underline{3}}^\alpha(q)
      + q^{n_1}\left( \chi_{4}^{\alpha-3}(q)+\chi_{\underline{4}}^{\alpha-3}(q)-\chi_{\underline{3}}^{\alpha-3}(q)   \right),\\
\chi_{4}^\alpha(q) \hspace{-1.7ex} & = & \hspace{-1.7ex} \chi_{\underline{3}}^\alpha(q) + q^{n_1} \chi_{\underline{4}}^{\alpha-4}(q),\\
\chi_{\underline{4}}^\alpha(q) \hspace{-1.7ex} & = & \hspace{-1.7ex} \chi_{\underline{3}}^\alpha(q) + q^{n_1}
      \chi_{4}^{\alpha-\underline{4}}(q),\\
\chi_{\underline{3}}^\alpha(q) \hspace{-1.7ex} & = & \hspace{-1.7ex} \chi_{\underline{2}}^\alpha(q) + q^{n_1}
      \chi_{\underline{2}}^{\alpha-\underline{3}}(q),\\
\chi_{\underline{2}}^\alpha(q) \hspace{-1.7ex} & = & \hspace{-1.7ex} \chi_{0}^\alpha(q) + q^{n_1}
      \chi_{0}^{\alpha-\underline{2}}(q).
\end{eqnarray*}
\end{prop}
The proof is similar to the proof in the $A_\ell$-case.

Set
\begin{eqnarray*}
\baza_{\Gamma'';0} \hspace{-1.7ex} & = & \hspace{-1.7ex} \{x_{\gamma_n}(-r_n)\cdots x_{\gamma_1}(-r_1) \in \baza_{\Gamma''} \st r_1\geq 2\},\\
\baza_{\Gamma'';\underline{3}} \hspace{-1.7ex} & = & \hspace{-1.7ex} \{x_{\gamma_n}(-r_n)\cdots x_{\gamma_1}(-r_1) \in \baza_{\Gamma''}
\st r_1\geq 2 \ \textrm{or}\ r_1=1, \gamma_1=\underline{3} \},\\
\baza_{\Gamma';0} \hspace{-1.7ex} & = & \hspace{-1.7ex} \{x_{\gamma_n}(-r_n)\cdots x_{\gamma_1}(-r_1) \in \baza_{\Gamma'} \st r_1\geq 2\},\\
\baza_{\Gamma';\underline{2}} \hspace{-1.7ex} & = & \hspace{-1.7ex} \{x_{\gamma_n}(-r_n)\cdots x_{\gamma_1}(-r_1) \in \baza_{\Gamma'}
\st r_1\geq 2 \ \textrm{or}\ r_1=1, \gamma_1=\underline{2} \},\\
\baza_{\Gamma';4,\underline{2}} \hspace{-1.7ex} & = & \hspace{-1.7ex} \{x_{\gamma_n}(-r_n)\cdots x_{\gamma_1}(-r_1) \in \baza_{\Gamma'}
\st r_1\geq 2 \ \textrm{or}\ r_1=1, \gamma_1\in\{4,\underline{2}\}\},\\
\baza_{\Gamma';\underline{4},\underline{2}} \hspace{-1.7ex} & = & \hspace{-1.7ex} \{x_{\gamma_n}(-r_n)\cdots x_{\gamma_1}(-r_1) \in \baza_{\Gamma'}
\st r_1\geq 2 \ \textrm{or}\ r_1=1, \gamma_1\in\{\underline{4},\underline{2}\}\},\\
\baza_{\Gamma';4,\underline{4},\underline{2}} \hspace{-1.7ex} & = & \hspace{-1.7ex} \{x_{\gamma_n}(-r_n)\cdots x_{\gamma_1}(-r_1) \in \baza_{\Gamma'}
\st r_1\geq 2 \ \textrm{or}\ r_1=1, \gamma_1\in\{4,\underline{4},\underline{2}\}\},
\end{eqnarray*}
and define $\baza_{\Gamma'';0}^\alpha$, $\baza_{\Gamma'';\underline{3}}^\alpha$, $\baza_{\Gamma';0}^\alpha$,
$\baza_{\Gamma';\underline{2}}^\alpha $, $\baza_{\Gamma';4,\underline{2}}^\alpha $, $\baza_{\Gamma';\underline{4},\underline{2}}^\alpha $,
$\baza_{\Gamma';4,\underline{4},\underline{2}}^\alpha $, and $\chi_{\Gamma'';0}^\alpha$,
$\chi_{\Gamma'';\underline{3}}^\alpha$, $\chi_{\Gamma';0}^\alpha$, $\chi_{\Gamma';\underline{2}}^\alpha$,
$\chi_{\Gamma';4,\underline{2}}^\alpha$, $\chi_{\Gamma';\underline{4},\underline{2}}^\alpha$,
$\chi_{\Gamma';4,\underline{4},\underline{2}}^\alpha$ in the obvious way.

Character formulas for these sets can be obtained in the same way as we did  for $\chi_{\Gamma'}^\alpha$
and $\chi_{\Gamma''}^\alpha$ in the previous section, by using character formulas for cases $A_3$, with $\omega=\omega_2$,
and $A_2$, with $\omega=\omega_2$, from Section \ref{ARR_sect}. We get
\begin{eqnarray}
\label{PrviKar_eq}\chi_{\Gamma'';0}^{\alpha''} \hspace{-1.7ex} & = & \hspace{-1.7ex} \frac{q^{n_1^2+n_4^2-n_1 n_4+n_1}}{(q)_{n_4}(q)_{n_1-n_4}},\\
\chi_{\Gamma'';\underline{3}}^{\alpha''} \hspace{-1.7ex} & = & \hspace{-1.7ex} \frac{q^{n_1^2+n_4^2-n_1 n_4+n_4}}{(q)_{n_4}(q)_{n_1-n_4}},\\
\chi_{\Gamma';0}^{\alpha'} \hspace{-1.7ex} & = & \hspace{-1.7ex} \frac{q^{n_4^2+n_1^2+n_3^2-n_4 n_1-n_1 n_3 +n_1}}
{(q)_{n_3} (q)_{n_1-n_3} (q)_{n_1-n_4} (q)_{n_4}},\\
\chi_{\Gamma';\underline{2}}^{\alpha'} \hspace{-1.7ex} & = & \hspace{-1.7ex} \frac{q^{n_4^2+n_1^2+n_3^2-n_4 n_1-n_1 n_3 }}
{(q)_{n_3} (q)_{n_1-n_3} (q)_{n_1-n_4} (q)_{n_4}} \left(q^{n_4+n_3}+q^{n_1} \frac{(1-q^{n_4}) (1-q^{n_3})}{1-q^{n_1}} \right),\\
\chi_{\Gamma';4,\underline{2}}^{\alpha'} \hspace{-1.7ex} & = & \hspace{-1.7ex} \frac{q^{n_4^2+n_1^2+n_3^2-n_4 n_1-n_1 n_3+n_3}}
{(q)_{n_3} (q)_{n_1-n_3} (q)_{n_1-n_4} (q)_{n_4}},\\
\chi_{\Gamma';\underline{4},\underline{2}}^{\alpha'} \hspace{-1.7ex} & = & \hspace{-1.7ex} \frac{q^{n_4^2+n_1^2+n_3^2-n_4 n_1-n_1 n_3+n_4}}
{(q)_{n_3} (q)_{n_1-n_3} (q)_{n_1-n_4} (q)_{n_4}},\\
\qquad\chi_{\Gamma';4,\underline{4},\underline{2}}^{\alpha'} \hspace{-1.7ex} & = & \hspace{-1.7ex}
\chi_{\Gamma';4,\underline{2}}^{\alpha'} + \chi_{\Gamma';\underline{4},\underline{2}}^{\alpha'} - \chi_{\Gamma';\underline{2}}^{\alpha'} \\
& = & \hspace{-1.7ex} \frac{q^{n_4^2+n_1^2+n_3^2-n_4 n_1-n_1 n_3}}{(q)_{n_3} (q)_{n_1-n_3} (q)_{n_1-n_4} (q)_{n_4}}
\left(1 - \frac{(1-q^{n_4}) (1-q^{n_3})}{1-q^{n_1}}\right).\label{ZadnjiKar_eq}
\end{eqnarray}
for $\alpha'$ and $\alpha''$ satisfying \eqref{alph1_jed} and \eqref{alph2_jed}, respectively.

\begin{prop} \label{D4PocUvDekomp_Prop}
Let $\underline{x},\underline{x}_1,\underline{x}_2$ be like in Proposition \ref{KarKomb_prop}. Then:
\begin{itemize}
\item[] $\underline{x}\in \baza_2 \Leftrightarrow  \underline{x}_1 \in \baza_{\Gamma'},  \underline{x}_2 \in \baza_{\Gamma''}$,
\item[] $\underline{x}\in \baza_3 \Leftrightarrow  \underline{x}_1 \in \baza_{\Gamma';4,\underline{4},\underline{2}},
\underline{x}_2 \in \baza_{\Gamma''}$,
\item[] $\underline{x}\in \baza_4 \Leftrightarrow  \underline{x}_1 \in \baza_{\Gamma';4,\underline{2}},  \underline{x}_2 \in \baza_{\Gamma''}$,
\item[] $\underline{x}\in \baza_{\underline{4}} \Leftrightarrow  \underline{x}_1 \in \baza_{\Gamma';\underline{4},\underline{2}},
\underline{x}_2 \in \baza_{\Gamma'';\underline{3}}$,
\item[] $\underline{x}\in \baza_{\underline{3}} \Leftrightarrow  \underline{x}_1 \in \baza_{\Gamma';\underline{2}},
 \underline{x}_2 \in \baza_{\Gamma'';\underline{3}}$,
\item[] $\underline{x}\in \baza_{\underline{2}} \Leftrightarrow  \underline{x}_1 \in \baza_{\Gamma';\underline{2}},
\underline{x}_2 \in \baza_{\Gamma'';0}$,
\item[] $\underline{x}\in \baza_0 \Leftrightarrow  \underline{x}_1 \in \baza_{\Gamma';0},  \underline{x}_2 \in \baza_{\Gamma'';0}$.
\end{itemize}
\end{prop}
The proposition can be proved by arguments similar to the ones used in the proof of  Proposition \ref{KarKomb_prop}.

Fix $n_1,n_2,n_3,n_4\geq 0$ satisfying \eqref{AlfaUvijet2_jed}, and
set $\alpha=n_1 \alpha_1+ n_2 \alpha_2 + n_3 \alpha_3 + n_4
\alpha_4$. Define $m_{\underline{2}}, m_{\underline{3}},
m_{\underline{4}}, m_{\underline{0}}$ by \eqref{m2n_jed}, $m'$ and
$m''$ by \eqref{MCrta_jed}, and $\alpha_i'$ and $\alpha_i''$, for
$i=0,\dots,n_4-m'-m''$, by \eqref{AlfaI_jed}. Proposition
\ref{D4PocUvDekomp_Prop} enables us to compute characters by using
analogues of formula \eqref{KarKomb_jed} and formulas
\eqref{PrviKar_eq}--\eqref{ZadnjiKar_eq}.

\begin{tm}
\begin{eqnarray}
\label{KarD4_jed}
\chi_\gamma^\alpha \hspace{-1.7ex} & = & \hspace{-1.7ex} \sum_{i=0}^{n_4-m'-m''} d_\gamma(\alpha) q^{f_i(\alpha)} \frac{(q)_{n_1-n_2+n_3-n_4+2(i+m')}}{
(q)_{n_3-n_4+i+m'} (q)_{n_1-n_2+i+m'} (q)_{i+m'} } \\
  &  & \hspace{5ex} \cdot\frac{1}{(q)_{n_1-n_2+n_3-n_4+i+m'} (q)_{n_4-i-m'} (q)_{n_2-n_3-i-m'}    }, \nonumber
\end{eqnarray} where
$f_i(\alpha)$ is defined by \eqref{KvadD_jed}, and $d_\gamma(\alpha)$ is defined by
$$
d_\gamma(\alpha) = \left\{ \begin{array}{ll}
1, & \textrm{for}\ \gamma=2,\\
\left(1-\frac{(1-q^{n_3-n_4+i+m'})(1-q^{i+m'})}{1-q^{n_1-n_2+n_3-n_4+2i+2m'}} \right), & \textrm{for}\ \gamma=3,\\
q^{n_3}, & \textrm{for}\ \gamma=4,\\
q^{n_4}, & \textrm{for}\ \gamma=\underline{4},\\
\left(q^{n_3+i+m'}-q^{n_1-n_2+n_3+i+m'}\frac{(1-q^{n_3-n_4+i+m'})(1-q^{i+m'})}{1-q^{n_1-n_2+n_3-n_4+2i+2m'}} \right) , & \textrm{for}\ \gamma=\underline{3},\\
\left(q^{n_2}-q^{n_1}\frac{(1-q^{n_3-n_4+i+m'})(1-q^{i+m'})}{1-q^{n_1-n_2+n_3-n_4+2i+2m'}} \right) , & \textrm{for}\ \gamma=\underline{2},\\
q^{n_1}, & \textrm{for}\ \gamma=0.
\end{array}\right.
$$
\end{tm}


\begin{thebibliography}{}


\bibitem[B1]{B1} I. Baranovi\'c, \textit{Combinatorial bases of Feigin-Stoyanovsky's type subspaces of level 2
standard modules for $D_4^{(1)}$}, math.QA/0903.0739

\bibitem[B2]{B2} I. Baranovi\'c, in preparation

\bibitem[C1]{C1} C. Calinescu, \textit{Intertwining vertex operators and certain representations of $\widehat{\mathfrak{sl}(n)}$},
Commun. Contemp. Math. {\bf 10} (2008), 47--79.

\bibitem[C2]{C2} C. Calinescu, \textit{Principal subspaces of higher-level standard $\widehat{\mathfrak{sl}(3)}$-modules},
J. Pure Appl. Algebra {\bf 210} (2007), 559--575.

\bibitem[CalLM1]{CalLM1} C. Calinescu, J. Lepowsky, A. Milas, \textit{Vertex-algebraic structure of the principal subspaces of
certain $A_1^{(1)}$-modules, I: level one case}, Int. J. Math. {\bf 19} (2008), 71--92.

\bibitem[CalLM2]{CalLM2} C. Calinescu, J. Lepowsky, A. Milas, \textit{Vertex-algebraic structure of the principal subspaces of
certain $A_1^{(1)}$-modules, II: higher-level case}, J. Pure Appl. Algebra, {\bf 212} (2008), 1928--1950

\bibitem[CalLM3]{CalLM3} C. Calinescu, J. Lepowsky, A. Milas, \textit{Vertex-algebraic structure of the principal subspaces of
level one modules for the untwisted affine Lie algebras of types A,D,E}, to appear in Journal of Algebra, math.QA/0908.4054

\bibitem [CLM1]{CLM1}    S. Capparelli, J. Lepowsky and A. Milas,
\textit{The Rogers-Ramanujan recursion and intertwining operators},
Comm. Contemporary Math. {\bf 5} (2003), 947--966.

\bibitem [CLM2]{CLM2}    S. Capparelli, J. Lepowsky, A. Milas,
\textit{The Rogers-Selberg recursions, the Gordon-Andrews identities and intertwining operators},  Ramanujan J.  {\bf 12}  (2006),  no. 3, 379--397

\bibitem[FJLMM]{FJLMM} B. Feigin, M. Jimbo, S. Loktev, T. Miwa and E. Mukhin, \textit{Bosonic formulas
for $(k,\ell)$-admissible partitions},  Ramanujan J. {\bf 7}  (2003),  no. 4, 485--517.; \textit{Addendum to
`Bosonic formulas for $(k,\ell)$-admissible partitions'},  Ramanujan J. {\bf 7}  (2003),  no. 4, 519--530

\bibitem[FJMMT1]{FJMMT1} B. Feigin, M. Jimbo, T. Miwa, E. Mukhin and Y. Takeyama, \textit{Fermionic formulas
for $(k,3)$-admissible configurations}, Publ. RIMS {\bf 40} (2004), 125--162.

\bibitem[FJMMT2]{FJMMT2} B. Feigin, M. Jimbo, T. Miwa, E. Mukhin and Y. Takeyama, \textit{Particle content of
the $(k,3)$-configurations}, Publ. RIMS {\bf 40} (2004), 163--220.

\bibitem[FS]{FS} A. V. Stoyanovsky and B. L. Feigin, \textit{Functional models of the representations of
current algebras, and semi-infinite Schubert cells}, (Russian) Funktsional. Anal. i Prilozhen. {\bf 28} (1994), no. 1, 68--90, 96;
translation in Funct. Anal. Appl. {\bf 28} (1994), no. 1, 55--72; preprint B. Feigin and A. Stoyanovsky,
\textit{Quasi-particles models for the representations of Lie algebras and geometry of flag manifold}, hep-th/9308079, RIMS 942.

\bibitem [G]{G} G. Georgiev,
\textit{Combinatorial constructions of modules for
infinite-dimensional Lie algebras, I. Principal subspace},  J. Pure
Appl. Algebra \textbf{112} (1996), 247--286.

\bibitem[J1]{J1} M. Jerkovi\'c,
\textit{Recurrence relations for characters of affine
Lie algebra $A_{\ell}^{(1)}$}, J. Pure Appl. Algebra {\bf 213},
913--926.

\bibitem[J2]{J2} M. Jerkovi\'c, in preparation

\bibitem[J3]{J3} M. Jerkovi\'c, PhD thesis, University of Zagreb, 2007.

\bibitem[K]{kac}
V.G. Kac,
\textit{Infinite-dimensional Lie algebras}, 3rd ed. Cambridge
University Press, Cambridge, 1990.

\bibitem[P1]{p1} M. Primc, \textit{ Vertex operator construction of standard modules
for $A_n^{(1)}$}, Pacific J. Math {\bf 162} (1994), 143--187.

\bibitem[P2]{p2} M. Primc, \textit{ Basic Representations sor classical
affine Lie algebras}, J. Algebra {\bf 228} (2000), 1--50.

\bibitem[P3]{p3} M. Primc, \textit{ $(k,r)$-admissible configurations and intertwining
operators}, Contemp. Math. {\bf 442} (2007), 425--434.

\bibitem[T1]{T1} G. Trup\v cevi\' c, \textit{Combinatorial bases of Feigin-Stoyanovsky's type subspaces of
level $1$ standard $\tilde{\gsl}(\ell+1,\C)$-modules}, to appear in Comm. Algebra, math.QA/0807.3363

\bibitem[T2]{T2} G. Trup\v cevi\' c, \textit{Combinatorial bases of Feigin-Stoyanovsky's type subspaces of
higher-level standard $\tilde{\gsl}(\ell+1,\C)$-modules}, J. Algebra {\bf 322} (2009), 3744--3774


\end{thebibliography}
\end{document}